\documentclass[12pt,a4paper]{article}
\usepackage[utf8]{inputenc}
\usepackage{mathrsfs,amsthm,xcolor,verbatim,bbm,amsmath,amsfonts,amssymb,nicefrac,enumitem,hyperref,bm,mathtools,xparse,etoolbox,graphicx,subcaption,siunitx}
\usepackage[a4paper, left=2cm, right=2cm]{geometry}
\usepackage{algorithm,algpseudocode}
\usepackage[nocompress]{cite}

\usepackage{todonotes}

\usepackage[english]{babel}
\usepackage{float} 
\usepackage{tikz-cd} 
\usepackage{booktabs} 

\usepackage[capitalise,sort]{cleveref} 

\crefname{enumi}{item}{items}
\crefname{equation}{}{}
\crefname{figure}{Figure}{Figures}
\crefname{listing}{Source code}{Source codes}
\crefname{lstlisting}{Source code}{Source codes}
\crefname{cor}{Corollary}{Corollaries}
\crefname{definition}{Definition}{Definitions}

\crefname{subsection}{Subsection}{Subsections}

\hypersetup{
    colorlinks,
    linkcolor={blue!80!black},
    citecolor={green},
    urlcolor={blue!80!black}
}

\usepackage{tcolorbox}
\tcbuselibrary{skins, breakable, listings}  
\newcounter{algorithmcounter}  
\renewcommand{\thealgorithmcounter}{\arabic{algorithmcounter}}  
\newtcolorbox[use counter=algorithmcounter]{myalgorithm}[3][]{
	enhanced,
	breakable,
	fonttitle=\bfseries,
	title=Algorithm~\thealgorithmcounter: #2,
	label={#3},  
	label type=algorithmcounter,
	#1,
	colframe=black,
	colback=white,
	coltitle=black,
	colbacktitle=white,
	sharp corners,
	boxrule=0.5pt,
	boxsep=1mm,
	top=1mm,
	bottom=1mm,
	left=0mm,
	right=0mm
}

\newcommand{\mycomment}[1]{\hfill \textcolor{gray}{\textit{\# #1}}}

\crefname{algorithmcounter}{Algorithm}{Algorithms}
\Crefname{algorithmcounter}{Algorithm}{Algorithms}
\crefname{line}{line}{lines}
\Crefname{line}{Line}{Lines}

\newcommand{\defaultParamDim}{d}
\newcommand{\momentprocess}{{\bf m}}
\newcommand{\Momentprocess}{{\bf v}}
\newcommand{\numberofsteps}{N}

\newcommand{\nos}{n}

\newcommand{\avgPa}{\delta}

\newcommand{\defaultdataDim}{\mathscr{d}}

\newcommand{\defaultStochLoss}{\mathscr{l}}
\newcommand{\defaultStochGradient}{\mathscr{g}}

\hypersetup{
	colorlinks,
	linkcolor={blue!80!black},
	citecolor={green},
	urlcolor={blue!80!black}
}

\theoremstyle{plain}
\newtheorem{theorem}{Theorem}[section]

\newtheorem{prop}[theorem]{Proposition}

\newtheorem{definition}[theorem]{Definition}

\DeclareMathAlphabet{\mathpzc}{OT1}{pzc}{m}{it}

\DeclareFontEncoding{LS1}{}{}
\DeclareFontSubstitution{LS1}{stix}{m}{n}
\DeclareMathAlphabet{\mathscr}{LS1}{stixscr}{m}{n}

\newcommand{\E}{\mathbb{E}}
\renewcommand{\P}{\mathbb{P}}

\newcommand{\R}{\mathbb{R}}
\newcommand{\N}{\mathbb{N}}

\newcommand{\cF}{\mathcal{F}}

\DeclarePairedDelimiter{\norm}{\lVert}{\rVert}
\DeclarePairedDelimiter{\abs}{\lvert}{\rvert}
\DeclarePairedDelimiter{\rbr}{(}{)}
\DeclarePairedDelimiter{\br}{[}{]}
\DeclarePairedDelimiter{\cu}{\{}{\}}

\renewcommand{\d}{ \mathrm{d}}

\newcommand{\andq}{\text{and}\qquad}

\ExplSyntaxOn

\bool_new:N \g_noteobserve

\NewDocumentCommand{\setnote}{}{
	\bool_gset_true:N \g_noteobserve
}

\NewDocumentCommand{\setobserve}{}{
	\bool_gset_false:N \g_noteobserve
}

\NewDocumentCommand{\nobs}{ o }{
	\IfValueT{#1}{
		\str_if_eq:noTF {note} {#1} {
			\bool_gset_true:N \g_noteobserve
		} {
			\str_if_eq:noTF {Note} {#1} {
				\bool_gset_true:N \g_noteobserve
			} {
				\bool_gset_false:N \g_noteobserve
			}
		}
	}
	\bool_if:nTF { \g_noteobserve } {
		\bool_gset_false:N \g_noteobserve
		note
	} {
		\bool_gset_true:N \g_noteobserve
		observe
	}
	\IfValueF{#1}{~}
}

\NewDocumentCommand{\Nobs}{ o }{
	\IfValueT{#1}{
		\str_if_eq:noTF {note} {#1} {
			\bool_gset_true:N \g_noteobserve
		} {
			\str_if_eq:noTF {Note} {#1} {
				\bool_gset_true:N \g_noteobserve
			} {
				\bool_gset_false:N \g_noteobserve
			}
		}
	}
	\bool_if:nTF { \g_noteobserve } {
		\bool_gset_false:N \g_noteobserve
		Note
	} {
		\bool_gset_true:N \g_noteobserve
		Observe
	}
	\IfValueF{#1}{~}
}

\int_new:N \g_furthermore

\NewDocumentCommand{\Moreover}{ o o }{
	\IfValueT{#1}{
		\str_case:nn {#1} {
			{Furthermore} {\int_set:Nn {\g_furthermore} {0}}
			{Moreover} {\int_set:Nn {\g_furthermore} {1}}
			{In~addition} {\int_set:Nn {\g_furthermore} {2}}
			{note} {\bool_gset_true:N \g_noteobserve}
			{observe} {\bool_gset_false:N \g_noteobserve}
		}
		\IfValueT{#2}{
			\str_case:nn {#2} {
				{Furthermore} {\int_set:Nn {\g_furthermore} {0}}
				{Moreover} {\int_set:Nn {\g_furthermore} {1}}
				{In~addition} {\int_set:Nn {\g_furthermore} {2}}
				{note} {\bool_gset_true:N \g_noteobserve}
				{observe} {\bool_gset_false:N \g_noteobserve}
			}
		}
	}
	\int_case:nn { \int_mod:nn {\g_furthermore} {3} } {
		{ 0 } { Furthermore,~\nobs that}
		{ 1 } { Moreover,~\nobs that}
		{ 2 } { In~addition,~\nobs that}
	}
	\int_incr:N \g_furthermore
	\IfValueF{#1}{~}
}

\bool_new:N \g_hencetherefore

\NewDocumentCommand{\hence}{}{
	\bool_if:nTF { \g_hencetherefore } {
		\bool_gset_false:N \g_hencetherefore
		hence~
	} {
		\bool_gset_true:N \g_hencetherefore
		therefore~
	}
}

\NewDocumentCommand{\Hence}{}{
	\bool_if:nTF { \g_hencetherefore } {
		\bool_gset_false:N \g_hencetherefore
		Hence,~we~obtain~
	} {
		\bool_gset_true:N \g_hencetherefore
		Therefore,~we~obtain~
	}
}

\seq_new:N \g_cflist_loaded
\seq_new:N \g_cflist_pending

\NewDocumentCommand{\cfadd}{ m }
{
	\seq_if_in:NnF \g_cflist_loaded { #1 } {
		\seq_if_in:NnF \g_cflist_pending { #1 } {
			\seq_gput_right:Nn \g_cflist_pending { #1 }
		}
	}
}

\NewDocumentCommand{\cfconsiderloaded}{ m }{
	\seq_gput_right:Nn \g_cflist_loaded {#1}
}

\NewDocumentCommand{\cfremove}{ m }
{
	\seq_gremove_all:Nn \g_cflist_pending { #1 }
}

\NewDocumentCommand{\cfload}{ o }
{
	\seq_if_empty:NTF \g_cflist_pending {\unskip} {
		(cf.\ \cref{\seq_use:Nn \g_cflist_pending {,}})\IfValueTF{#1}{#1~}{\unskip}
		\seq_gconcat:NNN \g_cflist_loaded \g_cflist_loaded \g_cflist_pending
		\seq_gclear:N \g_cflist_pending
	}
}

\NewDocumentCommand{\cfclear} {} {
	\seq_gclear:N \g_cflist_loaded
	\seq_gclear:N \g_cflist_pending
}

\NewDocumentCommand{\cfout}{ o }
{
	\seq_if_empty:NTF \g_cflist_pending {\unskip} {
		(cf.\ \cref{\seq_use:Nn \g_cflist_pending {,}})\IfValueTF{#1}{#1~}{\unskip}
		\seq_gclear:N \g_cflist_pending
	}
}

\NewDocumentCommand{\ifnocf} { m } {
	\seq_if_empty:NT \g_cflist_pending { #1 }
}

\ExplSyntaxOff

\NewDocumentEnvironment{cproof}{m}
{\begin{proof}[Proof of \cref{#1}]}%
	{\noindent The proof of \cref{#1} is thus complete.
\end{proof}}

\NewDocumentEnvironment{cproof2}{m}
{\begin{proof}[Proof of \cref{#1}]}%
	{\noindent This completes the proof of \cref{#1}.
\end{proof}}

\hypersetup{colorlinks=true}

\makeatletter
\@namedef{subjclassname@2020}{%
	\textup{2020} Mathematics Subject Classification}
\makeatother

\ExplSyntaxOn

\bool_new:N \g_forexample

\NewDocumentCommand{\eg}{ o }{
\IfValueT{#1}{
\str_if_eq:noTF {fe} {#1} {
\bool_gset_true:N \g_forexample
} {\bool_gset_false:N \g_forexample}
}
\bool_if:nTF { \g_forexample } {
\bool_gset_false:N \g_forexample
for~example
}{
\bool_gset_true:N \g_forexample
for~instance
}
}

\ExplSyntaxOff

\makeatletter
\ExplSyntaxOn
\seq_new:N \g_abbrs
\prop_new:N \g_abbr_counts
\tl_new:N \l_abbr_count_tl

\NewDocumentCommand{\abbr}{m m O{#1} m m O{#4}}{
	\expandafter\newcommand\csname#3\endcsname[1][]{
		\seq_if_in:NnTF \g_abbrs {#1} {
			\prop_get:NnN \g_abbr_counts {#1} \l_abbr_count_tl
			\prop_gput:Nnx \g_abbr_counts {#1} {\int_eval:n {\l_abbr_count_tl + 1}}
			\hyperref[#1]{#1}
		} {
			\seq_gput_left:Nn \g_abbrs {#1}
			\prop_gput:Nnn \g_abbr_counts {#1} {1}
			\expandafter\gdef\csname#1@def\endcsname{#2}
			\phantomsection\label{#1}
			\str_if_eq:nnTF{##1}{}{\emph{#2}}{##1}~(\hyperref[#1]{#1})
		}
	}
	\expandafter\newcommand\csname#6\endcsname[1][]{
		\seq_if_in:NnTF \g_abbrs {#1} {
			\prop_get:NnN \g_abbr_counts {#1} \l_abbr_count_tl
			\prop_gput:Nnx \g_abbr_counts {#1} {\int_eval:n {\l_abbr_count_tl + 1}}
			\hyperref[#1]{#4}
		} {
			\expandafter\gdef\csname#1@def\endcsname{#5}
			\seq_gput_left:Nn \g_abbrs {#1}
			\prop_gput:Nnn \g_abbr_counts {#1} {1}
			\phantomsection\label{#1}
			\str_if_eq:nnTF{##1}{}{\emph{#5}}{##1}~(\hyperref[#1]{#4})
		}
	}
}

\ExplSyntaxOff
\makeatother

\abbr{ANN}{artificial neural network}{ANNs}{artificial neural networks}
\abbr{DOS}{deep optimal stopping}{DOSs}{deep optimal stoppings}
\abbr{DL}{deep learning}{DLs}{deep learnings}
\abbr{iid}{independent and identically distributed}{iids}{deep learnings}
\abbr{SGD}{stochastic gradient descent}{SGDs}{stochastic gradient descent}
\abbr{Adam}{adaptive moment estimation}{Adams}{stochastic gradient descent}
\abbr{Padam}{parallel averaged \Adam}{Padams}{stochastic gradient descent}
\abbr{AdamW}{\Adam\ with weight decay}{AdamWs}{\Adams with weight decay}
\abbr{EMA}{exponential moving average}{EMAs}{exponential moving averages}
\abbr{RPA}{Ruppert--Polyak averaging}{RPAs}{Ruppert--Polyak averagings}
\abbr{SDE}{stochastic differential equation}{SDEs}{stochastic differential equations}
\abbr{ODE}{ordinary differential equation}{ODEs}{ordinary differential equations}
\abbr{PADAM}{parallel averaged \Adam}{PADAMs}{partially arithmetically averaged \Adam s}
\abbr{PDE}{partial differential equation}{PDEs}{partial differential equations}
\abbr{GD}{gradient descent}{GDs}{gradient descent}
\abbr{DKM}{deep Kolmogorov method}{DKMs}{deep Kolmogorov methods}
\abbr{DK}{deep Kolmogorov}{DKs}{deep Kolmogorovs}
\abbr{AI}{artificial intelligence}{AIs}{artificial intelligences}
\abbr{PINN}{physics-informed neural network}{PINNs}{physics-informed neural networks}
\abbr{DG}{deep Galerkin}{DGs}{deep Galerkins}
\abbr{deep BSDE}{deep backward stochastic differential equation}[deepBSDE]{deep BSDEs}{deep backward stochastic differential equations}[deepBSDEs]
\abbr{DRM}{deep Ritz method}{DRMs}{deep Ritz methods}
\abbr{GELU}{Gaussian error linear unit}{GELUs}{Gaussian error linear units}
\abbr{ReLU}{rectified linear unit}{ReLUs}{rectified linear units}
\abbr{ResNet}{residual neural network}{ResNets}{residual neural networks}
\abbr{OC}{optimal control}{OCs}{optimal controls}
\abbr{OS}{optimal stopping}{OSs}{optimal stoppings}
\abbr{HJB}{Hamiltonian--Jacobi--Bellman}{HJBs}{Hamiltonian--Jacobi--Bellmans}
\abbr{SOP}{stochastic optimization problem}{SOPs}{stochastic optimization problems}

\title{PADAM: Parallel averaged Adam reduces the error for\\ stochastic optimization in scientific machine learning}

\author{Arnulf Jentzen$^{1,2}$,
Julian Kranz$^{3,4}$, and Adrian Riekert$^{5}$\bigskip\\
\small{$^1$ School of Data Science and Shenzhen Research Institute of Big Data,}\vspace{-0.1cm}\\
\small{The Chinese University of Hong Kong, Shenzhen (CUHK-Shenzhen),}\vspace{-0.1cm}\\
\small{China; e-mail: \texttt{ajentzen}\textcircled{\texttt{a}}\texttt{cuhk.edu.cn}}\smallskip\\
\small{$^2$ Applied Mathematics: Institute for Analysis and Numerics,}\vspace{-0.1cm}\\
\small{University of M\"unster, Germany; e-mail: \texttt{ajentzen}\textcircled{\texttt{a}}\texttt{uni-muenster.de}}\smallskip\\
\small{$^3$ Applied Mathematics: Institute for Analysis and Numerics,}\vspace{-0.1cm}\\
\small{University of M\"unster, Germany; e-mail: \texttt{julian.kranz}\textcircled{\texttt{a}}\texttt{uni-muenster.de}}\smallskip\\
\small{$^4$ Machine Learning and Data Engineering: Department of Information Systems,}\vspace{-0.1cm}\\
\small{University of M\"unster, Germany; e-mail: \texttt{julian.kranz}\textcircled{\texttt{a}}\texttt{uni-muenster.de}}\smallskip\\
\small{$^5$ Applied Mathematics: Institute for Analysis and Numerics,}\vspace{-0.1cm}\\
\small{University of M\"unster, Germany; e-mail: \texttt{ariekert}\textcircled{\texttt{a}}\texttt{uni-muenster.de}}
}

\date{\today}

\begin{document}
	
\newcommand{\numproblems}{$13$ }
\maketitle
	
\begin{abstract}
Averaging techniques such as Ruppert--Polyak averaging and 
exponential movering averaging (EMA) are powerful approaches 
to accelerate optimization procedures of 
stochastic gradient descent (SGD) optimization methods
such as the popular ADAM optimizer. 
However, depending on the specific optimization problem under consideration, 
the type and the parameters for the averaging need to be adjusted 
to achieve the smallest optimization error. 
In this work we propose an averaging approach, 
which we refer to as parallel averaged ADAM (PADAM), 
in which we compute parallely different averaged variants 
of ADAM and during the training process dynamically select the variant 
with the smallest optimization error. 
A central feature of this PADAM approach is that this procedure
requires no more gradient evaluations than the usual ADAM optimizer
as each of the averaged trajectories relies on the same underlying ADAM trajectory
and thus on the same underlying gradients.
We test the proposed PADAM optimizer in \numproblems stochastic optimization and
deep neural network (DNN) learning problems
and compare its performance with known optimizers from the literature 
such as standard SGD, momentum SGD, Adam with and without EMA, and
ADAM with weight decay (ADAMW). 
In particular, we apply the compared optimizers  
to physics-informed neural network (PINN), deep Galerkin (DG), 
deep backward stochastic differential equation (deep BSDE) 
and deep Kolmogorov (DK) approximations for boundary value partial differential equation (PDE) problems 
(such as heat, Black--Scholes, Burgers, Allen--Cahn, and Hamiltonian--Jacobi--Bellman equations)
from scientific machine learning, as well as to DNN approximations for optimal control (OC) 
and optimal stopping (OS) problems. 
In nearly all 
of the considered numerical examples PADAM achieves, sometimes among others and sometimes exclusively, essentially the smallest optimization error.
This work thus strongly suggest to consider PADAM in the context of scientific machine learning problems
and also motivates further research for adaptive averaging procedures within the training of DNNs. 
The {\sc Python} source codes for each of the numerical experiments in this work can be found on {\sc GitHub} at 
\href{https://github.com/deeplearningmethods/padam}{https://github.com/deeplearningmethods/padam}. 
\end{abstract}
	
\tableofcontents

\section{Introduction}

\DL[\emph{Deep learning}]\ methods 
have not only revolutionized the state of the art 
of data driven \AI\ 
(cf., \eg, \cite{Brownlanguagemodel2020,Yihengsummary2023,
Rombachetal2022arXiv,Sahariaimagen2022_IMAGEN,
ramesh2021zeroshottexttoimagegeneration_DALLE}) 
but have 
also fundamentally changed 
the way how we solve scientific models 
such as \PDE, \OC, and inverse problems 
(cf., \eg, the overview articles 
\cite{beck2020overview,Weinan2021,Germainetal2021arXivOverviewArticle,MR4795589,MR4457972}).

\DL\ schemes usually consist of a class of deep \ANNs\ that are trained by \SGD\ optimization methods. 
Often not the standard \SGD\ method is the employed 
\SGD\ optimization method but instead more sophisticated accelerated or adaptive \SGD\ methods 
such as the \Adam\ \cite{KingmaBa2014} and the \AdamW\ \cite{LoshchilovHutter2017arXiv} optimizers 
are employed (cf., \eg, also \cite{Ruder2017overview,JentzenBookDeepLearning2023,bach2024learning} for overviews).

Moroever, averaging techniques such as \RPA\ \cite{Ruppert1988,MR1071220} (cf.\ also \cite{MR1167814}) 
and \EMA\ (cf., \eg, \cite{AhnCutkosky2024arXiv}) 
compose powerful approaches to accelerate optimization procedures of 
\SGD\ optimization methods. 
The classical \RPA\ approach seems to perform well 
for \SOPs\ in which 
the stochastic data in the \SOP\ is (nearly) \iid\ and in which the underlying optimizer 
is convergent to a (local) minimizer in the optimisation landscape 
(cf., \eg, \cite{AhnCutkosky2024arXiv}) 
but typically not in the situation of deep \ANN\ learning problems (see, \eg, \cref{sec:numerics} below).
Numerical simulations and theoretical investigations suggest 
that the reason for this poor behaviour of the \RPA\ approach 
in deep learning optimization is the issue that 
in the training of deep \ANNs\ we have that typically 
the gradient flow and the associated \SGD\ optimization methods
seem to apparently neither converge to a local or global minimizer 
nor a saddle point in the optimization landscape but instead 
fail to converge at all and diverge to infinity 
(cf., \eg, \cite{GallonJentzenLindner2022,MR4243432,Vardi2022,Lyu2020,Kranz2025}).
This topic is also closely related to the existence 
and non-existence, respectively, 
of minimizers in \ANN\ optimisation landscapes; cf., \eg,
\cite{JentzenRiekert2022_JML,MR4832358,MR4745362,GallonJentzenLindner2022,MR4243432}.

\EMA\ combined with \SGD\ optimization methods, in turn, seems 
to frequently accelerate the underlying optimization method 
for \SOPs\ in which the random variables describing the data in the \SOP\ are 
(nearly) \iid. In particular,
our preliminary work \cite{DereichJentzenRiekert2025arXiv}
suggests that \EMA\ combined 
with \Adam\ consistently reduces in several scientific machine learning cases
the approximation error for \PDE\ and \OC\ problems 
in which a huge number of essentially \iid\ training samples 
are available due to pseudo random number generators. 
We also refer, \eg, to \cite{morales-brotons2024exponential,SandlerZhmoginov2023,Defazioetal2024arXiv,ZhangChoromanskaLeCun2014arXiv,BusbridgeRamapuram2023EMA,GuoJin2023SWA,IzmailovPodoprikhin2018SWA,Athiwaratkun2018Average}
and the references therein for articles that propose and test
\SGD\ methods involving suitable averaging techniques
and we refer, \eg, to
\cite{AhnCutkosky2024arXiv,MR4184368,DereichMuller_Gronbach2019,
MR4580893,GadatPanloup2017arXiv,AhnMagakyanCutkosky2024arXiv,
MandtHoffmanBlei2017arXiv,Li2025}
and the references therein for works that mathematically study
averaged variants of \SGD\ methods.

However, depending on the specific \SOP\ under consideration, 
different types/parameters for the averaging result in quite different optimization errors 
and it remains an open question how to choose the specific averaging 
to achieve the smallest optimization error. 
In this work we propose an averaging approach, 
which we refer to as \Padam, 
in which we compute parallely different averaged variants 
of \Adam\ and during the training process dynamically select the variant 
with the approximately smallest optimization error.
A central feature of this approach is, on the one side, that 
this procedure requires, particularly for learning problems with large \ANNs, only minor additional computing time as each 
of the averaged trajectories relies on the same underlying \Adam\ trajectory and thus on the same underlying gradients 
but, on the other side, that this procedure accomplishes in many scientific machine learning based optimization problems 
smaller approximation errors than the most widely used optimizers
such as standard \SGD, momentum \SGD, \Adam\ with and without \EMA, and \AdamW.

In \cref{sec:numerics} below we test \Padam\ in \numproblems stochastic optimization and deep \ANN\ learning problems
and compare its performance with known optimizers
such as standard \SGD, momentum \SGD, \Adam\ with and without \EMA, and \AdamW. 
In particular, we apply the compared optimizers 
\begin{enumerate}[label=(\roman*)]
\item 
to polynomial regression problems
(see \cref{subsec:poly_reg}), 

\item 
to deep \ANN\ approximations for explicitly given high-dimensional target functions
(see \cref{subsec:supervised}),

\item 
to
\begin{itemize}
\item
\PINN,
\item
\DG,
\item
\deepBSDE, and
\item
\DK
\end{itemize}
approximations
for boundary value \PDE\ problems 
(such as heat, Black--Scholes, Burgers, Allen--Cahn, and \HJB\ equations)
from scientific machine learning, as well as 
\item 
to DNN approximations for \OC\ and \OS\ problems. 
\end{enumerate}
In nearly all 
of the considered numerical examples \Padam\ achieves, 
sometimes among others and sometimes exclusively, 
essentially the smallest optimization error, especially,
in the situation of scientific machine learning problems 
where a huge number of essentially \iid\ training samples are available due to pseudo random number generators.
Taking this into account, 
we strongly suggest to consider \Padam\ in the context of scientific machine learning problems.
This work also motivates further research for suitable adaptive averaging procedures within the training of deep \ANNs.
The {\sc Python} source codes for each of the numerical experiments in this work can be found on {\sc GitHub} at 
\href{https://github.com/deeplearningmethods/padam}{https://github.com/deeplearningmethods/padam}.

\subsubsection*{Structure of this article}

The remainder of this work is structured in the following way.
In \cref{sec:methods} we first recall the notion of the standard \Adam\ optimizer
and, based on this, we specify the proposed \Padam\ approach in detail.
In \cref{sec:numerics} we apply \Padam\ to
\numproblems different stochastic optimization and deep \ANN\ learning problems
and compare the obtained approximation errors
with those of optimization methods from the literature
such as standard \SGD, \Adam\ with and without \EMA, and \AdamW.
We close this paper with a short conclusion in \cref{sec:conclusion}.

\section{Parallel averaged Adam optimization}
\label{sec:methods}

\subsection{Standard Adam optimizer}

In \cref{ssub:padam} below we describe
the proposed \Padam\ approach. The formulation of \Padam\ is
based on the ``standard'' \Adam\ optimizer \cite{KingmaBa2014} and,
in view of this, we briefly recall within this subsection
the description of standard \Adam.
The precise form of
\cref{def:adam} comes from \cite[Definition~2.1]{DereichJentzenRiekert2025arXiv}.

\begin{definition}[Standard \Adam\ optimizer]
\label{def:adam}
Let 
$ \defaultParamDim, \defaultdataDim \in \N $, 
$
  ( \gamma_n )_{ n \in \N } 
  \allowbreak 
  \subseteq \R
$, 
$
  ( J_n )_{ n \in \N } \subseteq \N 
$,
$ 
  ( \alpha_n )_{ n \in \N } \subseteq [0,1)
$, 
$
  ( \beta_n )_{ n \in \N } \subseteq [0,1)
$, 
$
  \varepsilon \in (0,\infty) 
$, 
let $ ( \Omega, \mathcal{F}, \P ) $ be a probability space, 
for every $ n, j \in \N $ let
$
  X_{ n, j } \colon \Omega \to \R^{ \defaultdataDim }
$
be a random variable, 
let 
$
  \defaultStochLoss \colon \R^\defaultParamDim \times \R^{ \defaultdataDim } \to \R
$ 
be differentiable, 
let 
$
  \defaultStochGradient = (\defaultStochGradient_1,\ldots, \defaultStochGradient_\defaultParamDim)  \colon \R^\defaultParamDim \times \R^{ \defaultdataDim } 
   \to \R^{ \defaultParamDim }
$ 
satisfy for all 
$ \theta \in \R^{ \defaultParamDim } $,
$ x \in \R^{ \defaultdataDim } $
that
\begin{equation}
\label{eq:Adam_generalized_gradient}
  \defaultStochGradient(\theta,x) = \nabla_{ \theta } \defaultStochLoss( \theta, x )
  ,
\end{equation}
and let 
$ 
  \Theta = (\Theta^{(1)}, \ldots, \Theta^{(\defaultParamDim)}) \colon \N_0 \times \Omega \to \R^\defaultParamDim
$
be a function. 
Then we say that $ \Theta $ is the 
\Adam\ process 
for $ \defaultStochLoss $
with hyperparameters 
$ ( \alpha_n )_{ n \in \N } $, 
$ ( \beta_n )_{ n \in \N } $, 
$ ( \gamma_n )_{ n \in \N } $, 
$ \varepsilon \in (0,\infty) $, 
batch-sizes 
$ ( J_n )_{ n \in \N } $, 
initial value $ \Theta_0 $, and data $ ( X_{ n, j } )_{ (n,j) \in \N^2 } $ 
if and only if there exist  
$ 
  \momentprocess = (\momentprocess^{(1)},\ldots,\momentprocess^{(\defaultParamDim)}) \colon \N_0\times \Omega \to \R^\defaultParamDim
$
and 
$ 
  {\bf v} = ({\bf v}^{(1)},\ldots,{\bf v}^{(\defaultParamDim)}) \colon \N_0\times \Omega \to \R^\defaultParamDim
$ 
such that for all $ n \in \N $, $ i \in \{ 1, 2, \dots, \defaultParamDim \} $ 
it holds that 
\begin{equation}
  \momentprocess_0 = 0, \qquad 
  \momentprocess_n = \alpha_n \, \momentprocess_{n-1} + (1-\alpha_n)\br*{\frac{1}{J_n}\sum_{j = 1}^{J_n}
  \defaultStochGradient( \Theta_{n-1}, X_{n,j} )},
\end{equation}
\begin{equation}
  {\bf v}_0 = 0, 
  \qquad 
  {\bf v}_n^{(i)} = \beta_n\,{\bf v}_{n-1}^{(i)} + (1-\beta_n)\br*{\frac{1}{J_n} 
  \sum_{j = 1}^{J_n} \defaultStochGradient_{i}(\Theta_{n-1}, X_{n,j})}^2, 
\end{equation}
\begin{equation}
  \andq 
  \Theta_n^{(i)} = \Theta_{n-1}^{(i)} 
  - 
  \gamma_n 
  \,
  {\textstyle 
    \left[
      \varepsilon 
      + 
      \left[ 
        \frac{ 
          {\bf v}_n^{ (i) } 
        }{ 
          ( 1 - \prod_{ k = 1 }^n \beta_k ) 
        } 
      \right]^{ \nicefrac{1}{2} } 
    \right]^{ - 1 } 
  } 
  \br*{
    \frac{ \momentprocess_n^{ (i) } }{ ( 1 - \prod_{ k = 1 }^n \alpha_k ) }
  } .
\end{equation}
\end{definition}

Estimates for the optimization error of the \Adam\ optimizer can, \eg, be found in \cite{Barakat_2021_cvg,li2023convergenceadamrelaxedassumptions,DereichJentzen2024arXiv_Adam,ReddiKale2019,Defossez2022}
and the references therein.

\subsection{Parallel averaged Adam optimizer}
\label{ssub:padam}

Within this subsection we employ
\cref{def:adam} above
to formulate the proposed \Padam\ optimizer and its implementation
in \cref{def:averaged_adam_parallel} and \cref{alg:adam_avg_p} 
below.

\newcommand{\pparam}{\mathfrak{K}}
\begin{definition}[\Padam\ optimizer]
	\label{def:averaged_adam_parallel}
	Let 
	$ \defaultParamDim, \defaultdataDim, \pparam \in \N $,
	$
	  ( \gamma_n )_{ n \in \N }
	  \allowbreak
	  \subseteq \R
	$, 
	$
	  ( J_n )_{ n \in \N } \allowbreak
	  \subseteq \N
	$,
	$ 
	  ( \alpha_n )_{ n \in \N } \subseteq [0,1)
	$, 
	$
	  ( \beta_n )_{ n \in \N } \subseteq [0,1)
	$, 
	$
	  ( \avgPa_{ n, k } )_{ ( n, k ) \in \N \times \{ 1, 2, \dots, \pparam \} } \subseteq \R
	$, 
	$
	  \varepsilon \in (0,\infty)
	$, 
	let $ ( \Omega, \mathcal{F}, \P ) $ be a probability space, 
	for every $ n, j \in \N $ let
	$
	  X_{ n, j } \colon \Omega \to \R^{ \defaultdataDim }
	$
	be a random variable, 
	let 
	$
	  \defaultStochLoss \colon \R^\defaultParamDim \times \R^{ \defaultdataDim } \to \R
	$ 
	be differentiable, 
	let 
	$
	\defaultStochGradient = (\defaultStochGradient_1,\ldots, \defaultStochGradient_\defaultParamDim)  \colon \R^\defaultParamDim \times \R^{ \defaultdataDim } 
	\to \R^{ \defaultParamDim }
	$ 
	satisfy for all 
	$ \theta \in \R^{ \defaultParamDim } $,
	$ x \in \R^{ \defaultdataDim } $
	that
	\begin{equation}
		\label{eq:Adam_generalized_gradient3b}
		\defaultStochGradient(\theta,x) = \nabla_{ \theta } \defaultStochLoss( \theta, x )
		,
	\end{equation}
	and let 
	$ 
	  \Theta
	  \colon \N_0 \times \Omega \to \R^\defaultParamDim
	$
	be a function. 
	Then we say that 
	$ \Theta $ is the \Padam\ process
	for $ \defaultStochLoss $
	with hyperparameters
	$ ( \alpha_n )_{ n \in \N } $,
	$ ( \beta_n)_{ n \in \N } $,
	$ ( \gamma_n )_{ n \in \N } $,
	$ ( \avgPa_{ \nos,j } )_{ (\nos,k) \in \N \times \{ 1, 2, \dots, \pparam \} } $,
	$ \varepsilon \in (0,\infty) $,
	batch sizes $ ( J_n )_{ n \in \N } $,
	initial value $ \Theta_0 $,
	and
	data $ ( X_{ n, j })_{ (n,j) \in \N^2 } $
	if and only if there exist
	$ 
	  \vartheta^k \colon \N_0 \times \Omega \to \R^\defaultParamDim
	$,
	$ k \in \{ 0, 1, \dots, \pparam \} $,
	and
	$
	  K \colon \N \times \Omega \to \{ 0, 1, \dots, \pparam \}
	$
	such that 
	\begin{enumerate}[label=(\roman*)]
		\item 
		it holds that $ \vartheta^0 $ is the \Adam\ process
		for $ \defaultStochLoss $
		with hyperparameters
		$ ( \alpha_n )_{ n \in \N } $,
		$ ( \beta_n )_{ n \in \N } $,
		$ ( \gamma_n )_{ n \in \N } $,
		$ \varepsilon \in (0,\infty) $,
		batch sizes $ ( J_n )_{ n \in \N } $, 
		initial value $ \Theta_0 $,
		and data $ ( X_{ n, j } )_{ (n,j) \in \N^2 } $,
		\item 
		it holds for all $ k \in \{ 1, 2, \dots, \pparam \} $, $ n \in \N $ that
		\begin{equation}
		\label{eq:averaging_in_Padam}
		  \vartheta^k_0 = \Theta_0
		  \qquad
		  \text{and}
		  \qquad
			\vartheta^k_n
			= 
			\avgPa_{ n, k }
			\vartheta^k_{ n - 1 }
			+
			( 1 - \avgPa_{ n, k } )
			\Theta_n
			,
		\end{equation}
		and
		\item
		it holds for all $ n \in \N $ that
		$ \Theta_n =\vartheta^{ K_n }_n$
		and
		\begin{equation}
		  \sum_{ j = K_n J_n + 1 }^{ (K_n+1) J_n }
		  \defaultStochLoss( \Theta_n, X_{ n, j } )
	      =
	      \min_{ k \in \{ 1, 2, \dots, \pparam \} }
	      \left[
	        \sum_{ j = k J_n + 1 }^{ (k+1) J_n }
  	        \defaultStochLoss( \vartheta^k_n, X_{ n, j } )
		  \right]
		\end{equation}
	\end{enumerate}
\end{definition}

Note that the parameters
$
  ( \avgPa_{ n, k } )_{ (n,k) \in \N \times \{ 1, 2, \dots, \pparam \} } \subseteq \R
$
in \cref{def:averaged_adam_parallel}
correspond to the averaging weights
used in the possibly non-autonomous
different \EMA\ channels of \Adam\ in \cref{eq:averaging_in_Padam}.
In \cref{alg:adam_avg_p} below we now describe
the \Padam\ approach algorithmically.

%

\begin{myalgorithm}{\Padam\ }{alg:adam_avg_p}
	\begin{algorithmic}[1]
		\Statex\textbf{Setting:} The mathematical objects introduced in \cref{def:averaged_adam_parallel}
		\Statex\textbf{Input:} 
		$ \numberofsteps \in \N $
		\Statex\textbf{Output:} 
		\Padam\ process $ \Theta_{ \numberofsteps } \in \R^\defaultParamDim $
		after $ \numberofsteps $ steps
		
		\vspace{-2mm} 
		\noindent\hspace*{-8.1mm}\rule{\dimexpr\linewidth+9.24mm}{0.4pt}  
		\vspace{-3mm} 
		
		\State $ \vartheta \gets \Theta_0 $
		\For{$ k \in \cu{ 1, 2, \dots, \pparam } $}
		\State $ \theta_k \gets \Theta_0 $
		\EndFor
		\State $ \momentprocess \gets 0 $
		\State $ \Momentprocess \gets 0 $
		\For{$ \nos \in \cu{ 1, 2, \dots, \numberofsteps } $}
		\State $ g \gets ( J_{ \nos } )^{ - 1 } \sum_{ j = 1 }^{ J_n } \defaultStochGradient( \vartheta, X_{ \nos, j } ) $
		\State $ \momentprocess \gets \alpha_{ \nos } \momentprocess + ( 1 - \alpha_{ \nos } ) g $
		\State $ \Momentprocess \gets \beta_{ \nos } \momentprocess + ( 1 - \beta_{ \nos } ) g^{ \otimes 2 } $ 
		\mycomment{Square $ g^{ \otimes 2 } $ is understood componentwise}
		\State $ \hat{\momentprocess} \gets \momentprocess / ( 1 - \prod_{ k = 1 }^{ \nos } \alpha_k ) $
		\State $ \hat{\Momentprocess} \gets \Momentprocess / ( 1 - \prod_{ k = 1 }^{ \nos } \beta_k ) $
		\State $ \vartheta \gets \vartheta - \gamma_{ \nos } \hat{\momentprocess} / ( \hat{\Momentprocess}^{ \otimes (1/2) } + \varepsilon ) $ 
		\mycomment{Root $ \Momentprocess^{ \otimes ( 1 / 2 ) } $ is understood componentwise}
		\For{$ k \in \cu{ 1, 2, \dots, \pparam } $}
		\State $ \theta_k \gets \avgPa_{ \nos,j } \theta_k + ( 1 - \avgPa_{ \nos,k } ) \vartheta $ 
		\mycomment{Update averaged iterates}
		\EndFor
		\EndFor
		\noindent 
		\State $k_{*} \gets 1$
		\For{$ k \in \cu{ 1, 2, \dots, \pparam } $}
		\If{$\sum_{ j = kJ_N + 1 }^{ (k+1) J_N }
			\defaultStochLoss( \theta_k, X_{ N, j } ) < \sum_{ j = k_{*}J_N + 1 }^{ (k_{*}+1) J_N } \defaultStochLoss( \theta_{k_{*}}, X_{ N, j } )$} 
		\State $k_{*} \gets k$ \mycomment{Choose optimal $k\in \{1,\dotsc,\pparam\}$}
		\EndIf
		\EndFor
		\State \Return $ \theta_{k_{*}} $
	\end{algorithmic}
\end{myalgorithm}

In the numerical experiments in \cref{sec:numerics} below
we tested the cases $ \pparam = 3 $ and $ \pparam = 10 $
for \cref{alg:adam_avg_p} and
referred to these methods as PADAM3 and PADAM10, respectively.
For displaying the performance of the \Padam\ algorithms,
we fixed a threshhold $ n_T \in \{ 500, 5000 \} $ and
compute the test errors $\mathcal L_{\mathrm{test}}(\Theta_{n,j})$ for the different \emph{channels} $\Theta_{n,1}, \dotsc,\Theta_{n,k}$ whenever $n$ is divisible by $n_T$ and then plot the test error of the best performing channel for the next $n_T$ gradient steps. We chose $n_T=500$ for the \OC\ and \OS\ problems (see \cref{subsec:OC} and \cref{subsec:OS}) and $n_T=5000$ for all the other problems.
Below we list the averaging parameters for PADAM3.
Here, $N$ is the total number gradient steps for the underlying Adam optimizer.
\begin{enumerate}[label=(\roman*)]
	\item $\avgPa_{n,1} = 0.999$,
	\item $\avgPa_{n,2} = 1 - n^{ - 0.7 } $,
	\item $\avgPa_{n,3} = 1 - 0.1 \exp( \frac{ - 2 n \ln( 10 ) }{ N } ) $
\end{enumerate}
Below, we list the averaging parameters for PADAM10. Here, $N$ is the total number gradient steps for the underlying Adam optimizer.
\begin{enumerate}[label=(\roman*)]
	\item $\avgPa_{n,1}=0.99$, 
	\item $\avgPa_{n,2}=0.999$,
	\item $\avgPa_{n,3}= 1-n^{-0.6}$,
	\item $\avgPa_{n,4}=1-n^{-0.7}$,  
	\item $\avgPa_{n,5}=1-n^{-0.8}$,  
	\item $\avgPa_{n,6}=1-0.5n^{0.7}$,
	\item $\avgPa_{n,7}=1-0.1 \exp( \frac{ - 2 n \ln(10)}{N} ) $,
	\item $\avgPa_{n,8}=1-0.01 \exp( \frac{ - n \ln(10)}{N} ) $,
	\item $\avgPa_{n,9}=1-0.1 \exp( \frac{ - 3 n \ln( 10 )}{N} ) $,
	\item $\avgPa_{n,10}=1.-0.1 \exp( \frac{ - 5 n \ln( 10 )}{N} ) $.
\end{enumerate}
The hyperparameters for these channels were found by trial and error.

\section{Numerical experiments}
\label{sec:numerics}

\subsection{Polynomial regression}
\label{subsec:poly_reg}

\begin{figure}[H]\caption{Polynomial Regression Problem
	}\label{Polynomial Regression_[25]_174133504683714}
	\includegraphics[width=\linewidth]{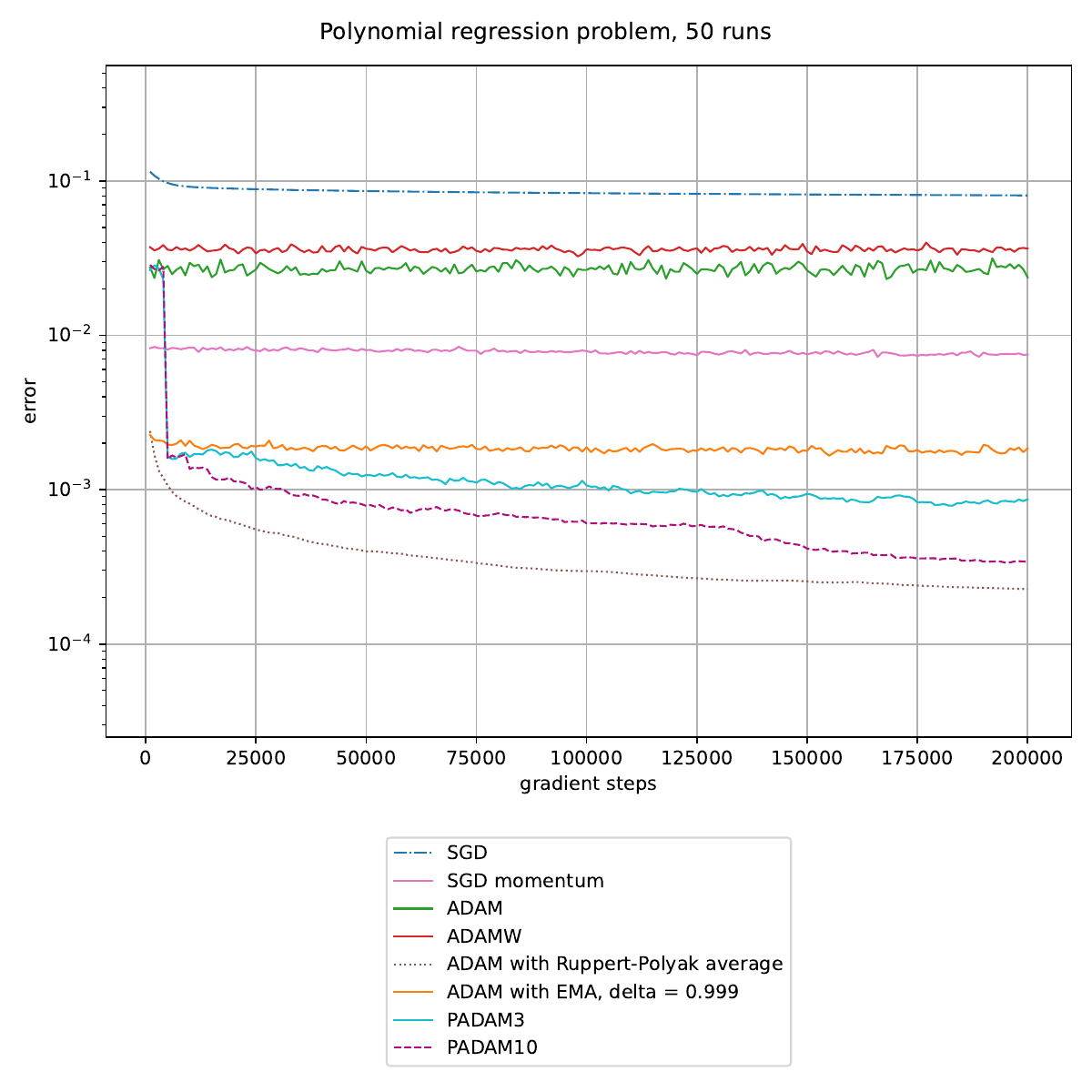}
\end{figure}

In our first numerical example
in Figure \ref{Polynomial Regression_[25]_174133504683714}
we consider the problem to approximate the explicitly given function
$ [-1,1] \ni x \mapsto \sin( \pi x ) \in \R $
in the $ L^2( [-1,1]; \R ) $-sense by means of polynomials of degree at most $ 25 $.
More formally, we aim to minimize the function
\begin{equation}
\textstyle
	\R^{d+1} \ni 
	\theta = ( \theta_0, \theta_1, \dots, \theta_d ) 
	\mapsto 
	\int\limits_{ - 1 }^1
	\Bigl(
	  \sin( \pi x ) - \sum\limits_{ k = 0 }^d \theta_k x^k
	\Bigr)^{ \! 2 } \, \d x
	\in \R
\end{equation}
for $ d = 25 $,
leading to an $ 26 $-dimensional convex optimization problem.
In the training we use mini-batches of size $ \num{256} $
and constant learning rates of size $ \num{0.01} $
and we add centered Gaussian noise with variance $ \num{0.2} $
to the output.
Furthermore,
in Figure \ref{Polynomial Regression_[25]_174133504683714}
we approximate the relative $ L^2( [-1,1]; \R ) $-error
through a Monte Carlo approximation with $ \num{50000} $ Monte Carlo samples
and we approximate the $ L^1 $-error with respect to the probability space
through a Monte Carlo approximation
with $ \num{50} $ independent simulations.

\subsection{Deep artificial neural network (ANN) approximations for Gaussian densities}
\label{subsec:supervised}

In the next example we consider the approximation
in the $ L^2( [-2,2]^d; \R ) $-sense
of the normal density function
\begin{equation}
	\label{eq:gauss}
	\R^d \ni x \mapsto
	\exp\bigl(
	  - \tfrac{ \| x \|^2 }{ 2 \sigma^2 }
	\bigr)
	\in \R
\end{equation}
in $ d = 20 $ dimensions with
the standard deviation parameter
$ \sigma = \sqrt{3} $.
In Figure \ref{supervised problem_[20, 300, 500, 100, 1]_174617865182029}
we approximate the function in \cref{eq:gauss}
using \ANNs\ with the \ReLU\ activation
(see, \eg, \cite[Subsection~1.2.3]{JentzenBookDeepLearning2023})
and three hidden layers consisting of $300$, $500$, and $100$ neurons, respectively.
For the input distribution we choose the continuous uniform distribution on $ [-2,2]^d $
and for the loss function we employ the mean squared error, both, for the training and the test loss.
In the training we employ mini-batches of size
$ \num{256} $ and
constant learning rates of size $ 10^{ - 4 } $.
Furthermore,
in Figure \ref{supervised problem_[20, 300, 500, 100, 1]_174617865182029}
we approximate the $ L^2( [-2,2]^d; \R ) $-error
through a Monte Carlo approximation with $ \num{100 000} $ Monte Carlo samples
and we approximate the $ L^1 $-error with respect to the probability space
through a Monte Carlo approximation
with $ 50 $ independent simulations.

\begin{figure}[H]\caption{\label{supervised problem_[20, 300, 500, 100, 1]_174617865182029}Supervised deep \ANN\ learning
of Gaussian densities}
\includegraphics[width=\linewidth]{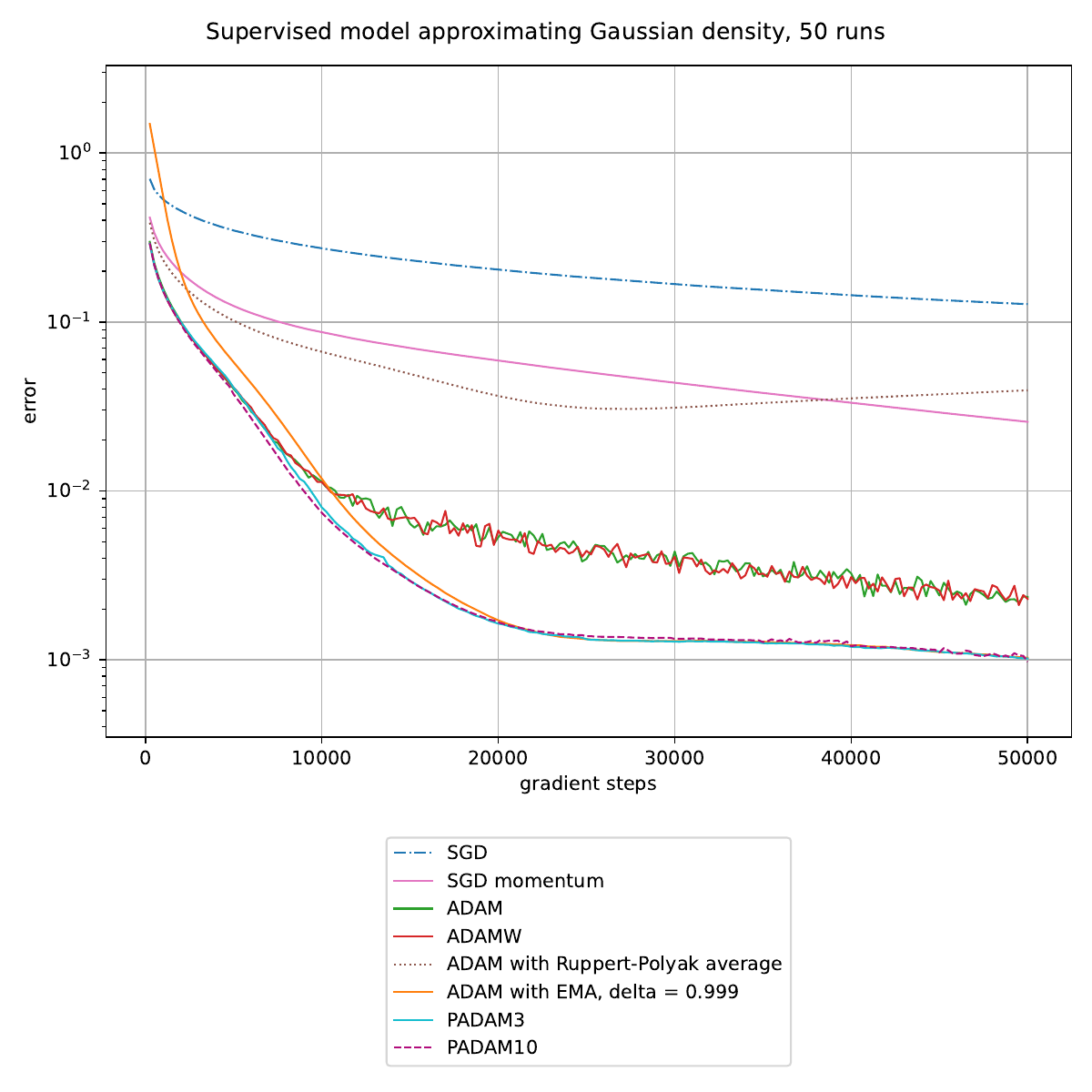}
\end{figure}

\subsection{Deep Kolmogorov method (DKM) for heat equation}
\label{subsec:DKM1}

In the next example we employ the \DKM\ from Beck et al.~\cite{Beck2021Kolmogorov}
to approximately solve the heat \PDE\ on $ \R^d $ for $ d = 10 $.
More formally, we aim to approximate the solution $ u \colon [0, T] \times \R^d \to \R $ of the
initial value \PDE\ problem
\begin{equation}
\label{eq:heat}
  \tfrac{ \partial u }{ \partial t }
  = \Delta_x u
  ,
\qquad
  u( 0, x ) = \| x \|^2
\end{equation}
for $ t \in [0,T] $, $ x \in \R^d $
at the final time $ T = 2 $
on the set $ [ -1, 1]^d \subseteq \R^d $.
The \PDE\ can be reformulated as a \SOP\ (cf.~Beck et al.~\cite{Beck2021Kolmogorov}) and
thus \SGD\ methods such as \Adam\ are applicable to compute an approximate minimizer.
In Figure~\ref{Heat Model_[10, 50, 100, 50, 1]_174109636785599}
we approximate the solution of \cref{eq:heat}
by fully connected feedforward \ANNs\ with the \GELU\ activation
(see, \eg, \cite[Subsection~1.2.6]{JentzenBookDeepLearning2023})
and three hidden layers consisting of $ 50 $, $ 100 $, and $ 50 $ neurons, respectively.
In the training we employ mini-batches of size $ 256 $ and constant learning rates of size $ 10^{ - 4 } $.
To compute the error in Figure~\ref{Heat Model_[10, 50, 100, 50, 1]_174109636785599}
we employ the fact that the exact solution
$ u \colon [0,T] \times \R^d \to \R $ of \cref{eq:heat}
satisfies that for all $ t \in [0,T] $, $ x \in \R^d $
it holds that
$ u(t,x) = \|x\|^2 + 2 d t $.
Furthermore,
in Figure~\ref{Heat Model_[10, 50, 100, 50, 1]_174109636785599}
we approximate the relative $ L^2( [-1,1]^d; \R ) $-error
through a Monte Carlo approximation with $ 10^5 $ Monte Carlo samples
and we approximate the
$ L^1 $-error with respect to the probability space
through a Monte Carlo approximation
with $ 50 $ independent simulations.

\begin{figure}[H]\caption{\label{Heat Model_[10, 50, 100, 50, 1]_174109636785599}Deep Kolmogorov method for Heat equations}
\includegraphics[width=\linewidth]{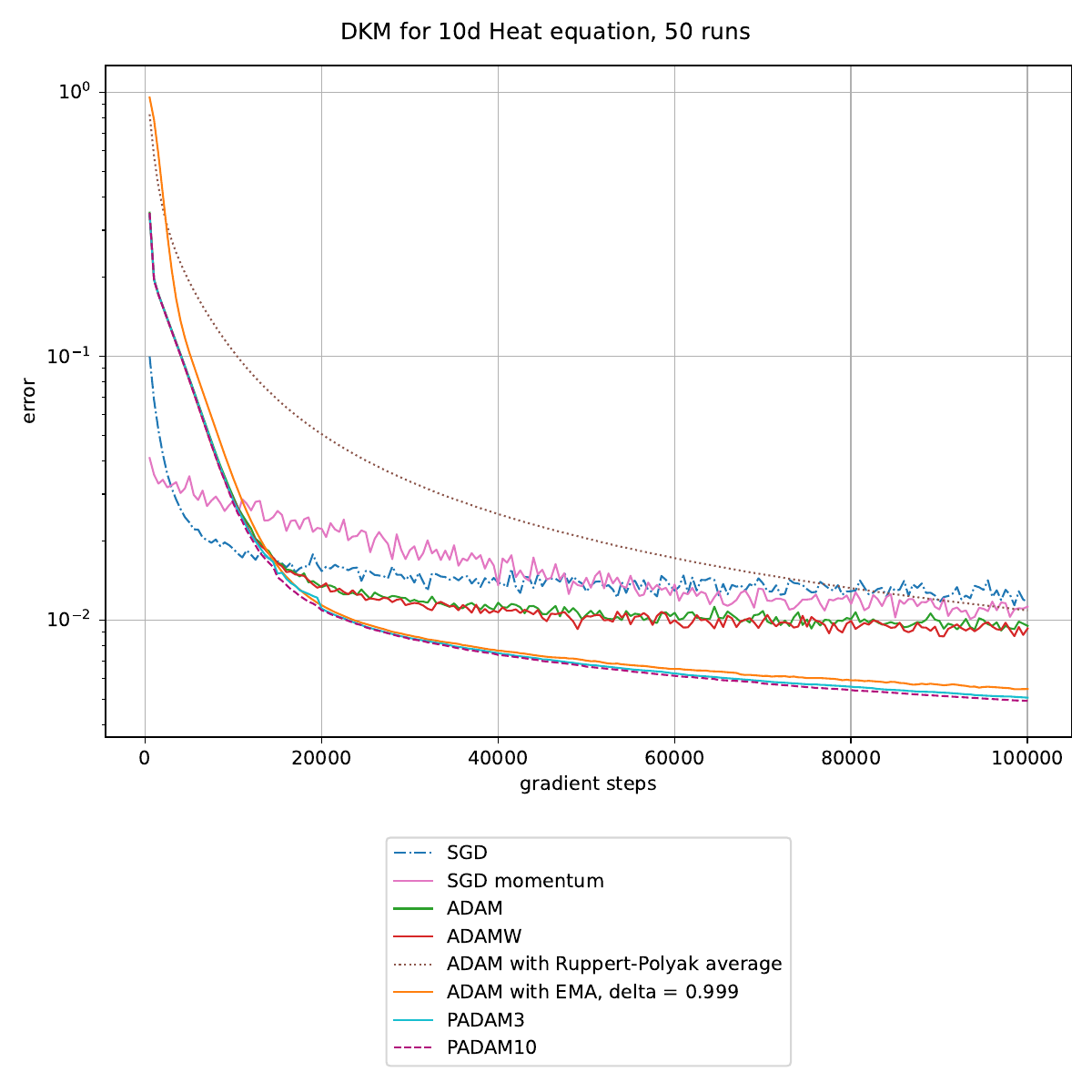}
\end{figure}

%
%

\subsection{DKM for Black--Scholes equation}
\label{subsec:DKM2}


In the next example we apply again the \DKM\ to approximately
compute the solution
$
  u \colon [0, T] \times \R^d \to \R
$
of the Black--Scholes \PDE
\begin{equation}
\label{eq:black-scholes}
\begin{split}
  \tfrac{\partial u }{\partial t }
&
\textstyle
  =
  \frac{ 1 }{ 2 }
  \left[
    \sum\limits_{ i = 1 }^d
    |
      { \sigma_i x_i }
    |^2
    \tfrac{ \partial^2 u }{ \partial x_i^2 }
  \right]
  +
  \mu
  \left[
    \sum\limits_{ i = 1 }^d
    x_i \tfrac{\partial u }{ \partial x_i }
  \right]
,
\\
  u ( 0 , x )
&
  = \exp( - r T )
  \max\{ \max\{ x_1, x_2, \dots, x_d \} - K , 0 \} ,
\end{split}
\end{equation}
for $ t \in [0,T] $, $ x = ( x_1, \dots, x_d ) \in \R^d $
at the final time $ T = 1 $
on $ [ 90, 110 ]^d \subseteq \R^d $
where $ d = 10 $,
where
$ \sigma_1, \sigma_2, \dots, \sigma_d \in \R $
satisfy for all $ i \in \{ 1, 2, \dots, d \} $ that
$ \sigma_i = \frac{ i + 1 }{ 2 d } $,
where
$ r = - \mu = \frac{ 1 }{ 20 } $,
and where
$ K = 100 $.
In Figure~\ref{Black Scholes Model_[10, 200, 300, 200, 1]_174133504617140}
we approximate the solution of \cref{eq:black-scholes}
by means of
fully connected feedforward \ANNs\ with the \GELU\ activation
and three hidden layers consisting of $ 200 $, $ 300 $, and $ 200 $ neurons, respectively,
and a batch normalization layer before the first hidden layer.
In the training we use mini-batches of size $ 256 $ and constant learning rates of size $ 10^{ - 4 } $.
To approximately compute the error in Figure~\ref{Black Scholes Model_[10, 200, 300, 200, 1]_174133504617140},
we use the Monte Carlo method with $ \num{10240000} $ Monte Carlo samples
based on the Feynman--Kac formula for \cref{eq:black-scholes} to approximate
the unknown exact solution $ u \colon [0,T] \times \R^d \to \R $ of \cref{eq:black-scholes}.
Furthermore,
in Figure~\ref{Black Scholes Model_[10, 200, 300, 200, 1]_174133504617140}
we approximate
the relative $ L^2( [90,110]^d; \R ) $-error
through a Monte Carlo approximation
with $ 10^5 $ Monte Carlo samples
and we approximate the $ L^1 $-error with respect
to the probability space through a Monte Carlo approximation
with $ 50 $ independent simulations.

\begin{figure}[H]\caption{\label{Black Scholes Model_[10, 200, 300, 200, 1]_174133504617140}Deep Kolmogorov method for Black Scholes equations}
\includegraphics[width=\linewidth]{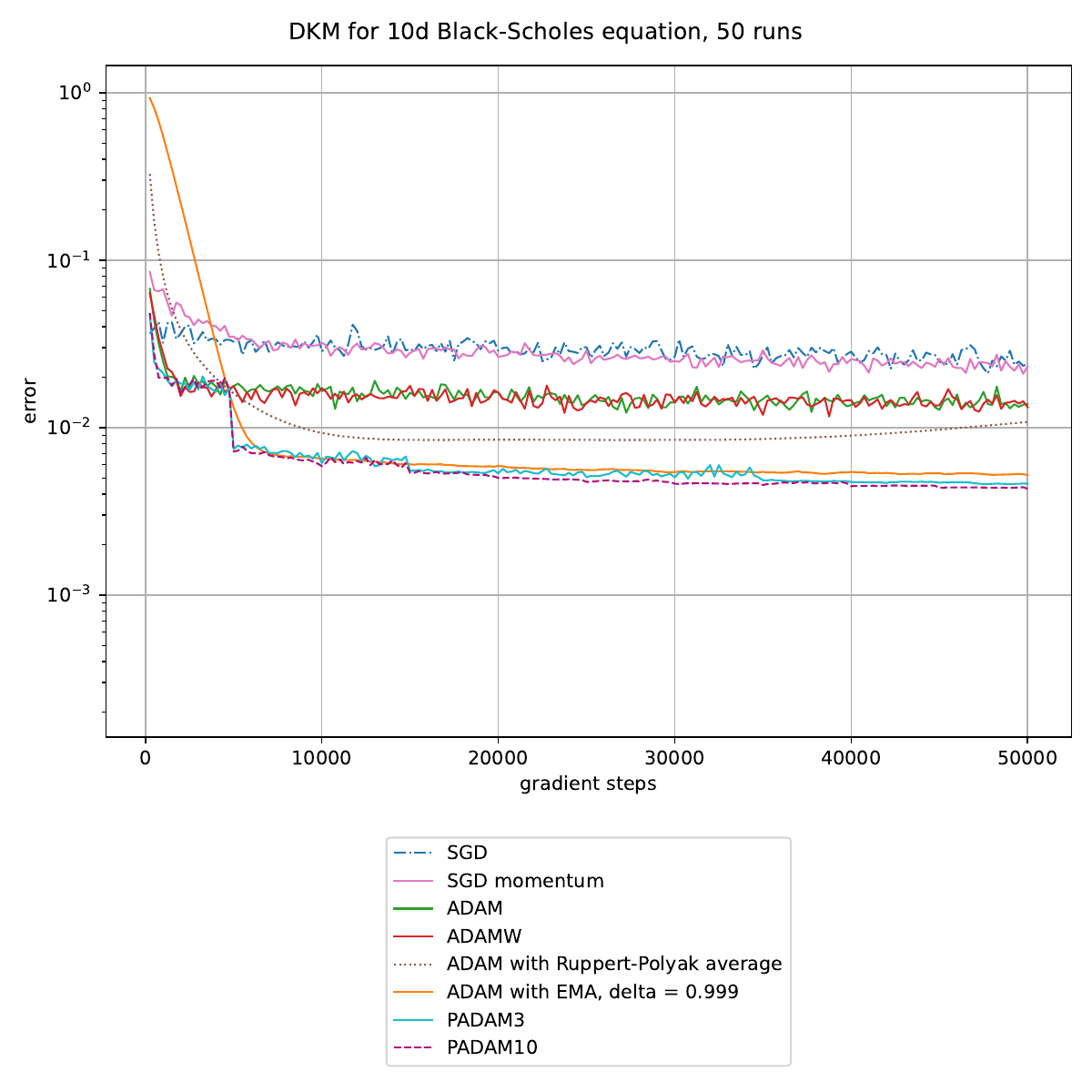}
\end{figure}
%
%

\subsection{Quadratic stochastic minimization problem}
\label{subsec:quadratic}

In the next example we aim to minimize the quadratic function
\begin{equation}
\label{eq:quadratic}
  \R^d \ni \theta \mapsto \E\bigl[ \norm{ \theta - X }^2 \bigr] \in \R
\end{equation}
where $ X $ is a $ d $-dimensional standard normal random variable
and where $ d = 10 $.
Note that \cref{eq:quadratic} is a strongly convex function
with the unique global minimum at $ \theta = \E[ X ] = 0 $.
Taking this into account,
in Figure~\ref{Quadratic Minimization Problem, d=10, var=1, mean=0_[10]_174220121440214}
we consider the error
$
  \R^d \ni \theta \mapsto \frac 1 d \|\theta \|^2 = \frac 1 d \| \theta - \E[X] \|^2 \in \R
$.
In the training we use mini-batches of size $ \num{256} $.
Moreover, for \SGD\ and \SGD\ with momentum we employ constant learning rates of size $ \num{0.001} $
and for all other optimization methods
we employ constant learning rates of size $ \num{0.01} $.
Furthermore, we approximate the $L^1$-error with respect to the probability space
through a Monte Carlo approximation
with $ 50 $ independent simulations.

\begin{figure}[H]\caption{Quadratic problem
}\label{Quadratic Minimization Problem, d=10, var=1, mean=0_[10]_174220121440214}
\includegraphics[width=\linewidth]{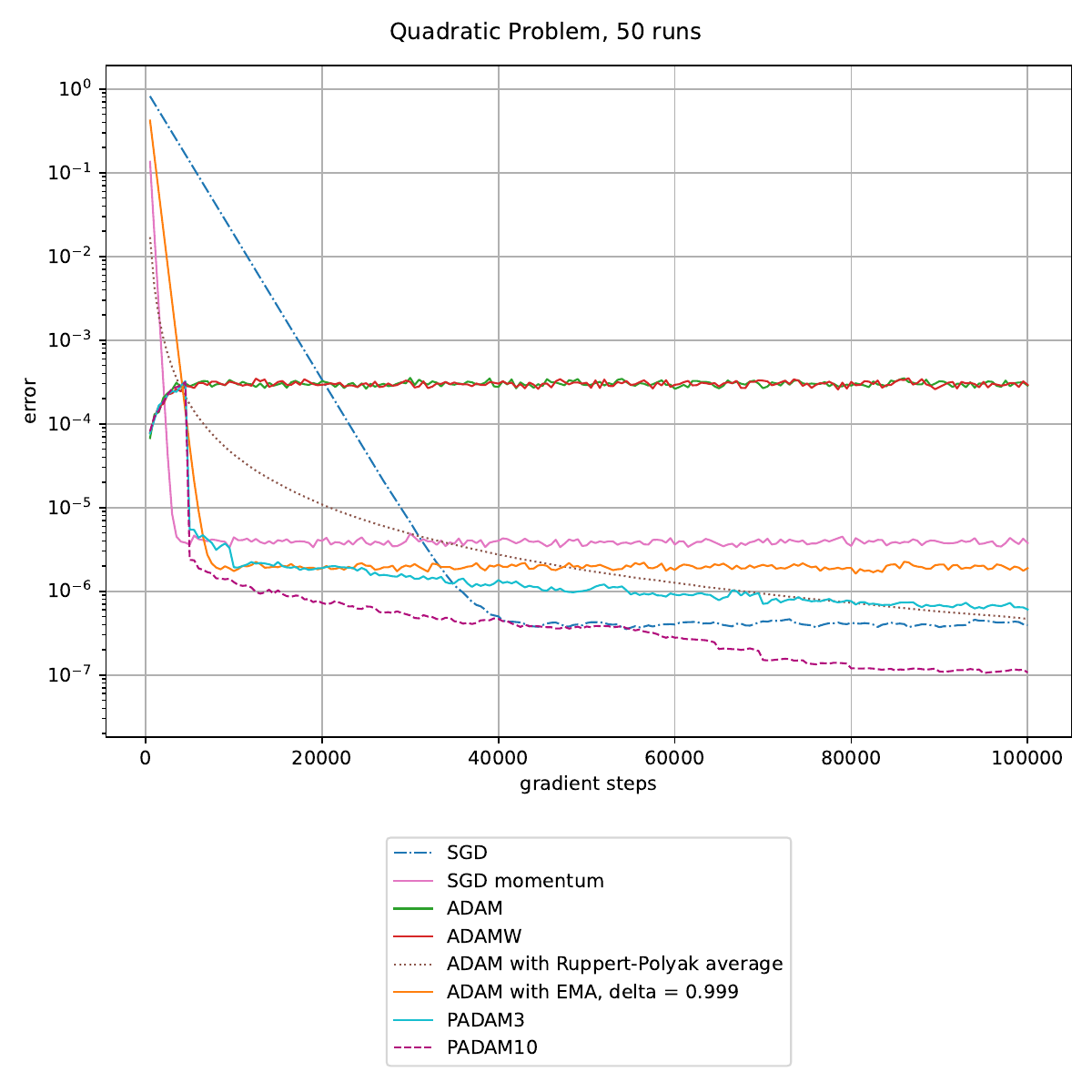}
\end{figure}

\subsection{Deep Ritz for Poisson equation}
\label{subsec:DRM1}

\begin{figure}[H]\caption{Deep Ritz for a Poisson equation
}\label{Deep Ritz Heat Model_None_174167699975033}
\includegraphics[width=\linewidth]{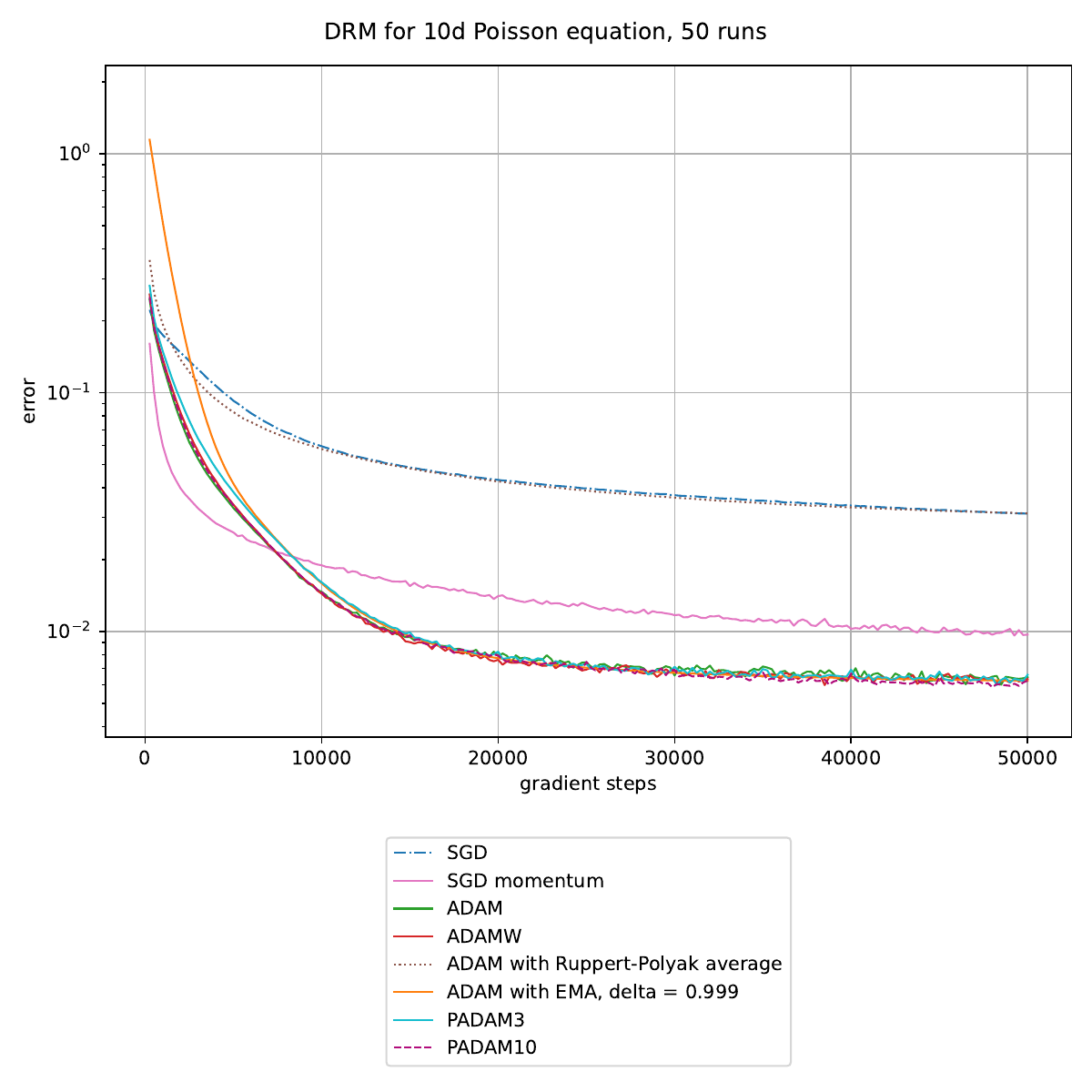}
\end{figure}

In the next example problem
we employ the \DRM\ (cf.~\cite{EYu2018}) to approximately compute
the solution $ u \colon [-1 , 1 ]^d \to \R $
of the $ d $-dimensional Poisson equation
\begin{equation}
\label{eq:Poisson}
\begin{split}
  \Delta u( x ) = 2 d ,
\qquad
  u(y) = \norm{y}^2
\end{split}
\end{equation}
for $ x \in ( -1, 1 )^d $, $ y \in \partial ( -1, 1 )^d $
where $ d = 10 $.
The exact solution of $ u \colon [-1 , 1 ]^d \to \R $
of \cref{eq:Poisson} satisfies that for all
$ x = ( x_1, \dots, x_d ) \in [-1,1]^d $ it holds that
$ u( x ) = \norm{x}^2 = \sum_{i=1}^d ( x_i )^2 $.
In Figure~\ref{Deep Ritz Heat Model_None_174167699975033}
we approximate the exact solution of \cref{eq:Poisson}
with fully connected feedforward \ANNs\ with the \GELU\ activation
and $ 6 $ hidden layers each consisting of $ 32 $ neurons.
The boundary condition is incorporated using a penalty method.
For the training we employ mini-batches of size $ 1024 $.
Moreover, for \SGD\ and \SGD\ with momentum we employ constant learning rates of size $ 3 \cdot 10^{ - 6 } $
(to avoid divergence) and for all other optimization methods in Figure~\ref{Deep Ritz Heat Model_None_174167699975033}
we employ constant learning rates of size $ 3 \cdot 10^{ - 4 } $.
In Figure~\ref{Deep Ritz Heat Model_None_174167699975033}
we approximate the relative $ L^2( [-1,1]^d; \R ) $-error
through a Monte Carlo approximation with $ 10^4 $ Monte Carlo samples
and we approximate the $ L^1 $-error with respect to the probability space
through a Monte Carlo approximation with $ 50 $ independent simulations.

\subsection{Deep Ritz for p-Laplace equation}
\label{subsec:DRM2}

In the next example, which is based on \cite[Section 5.3]{DondlMuller2022},
we apply the \DRM\ to
the nonlinear \PDE
\begin{equation}
\label{eq:nonlinear_DRM}
  - \operatorname{div} \rbr*{ \norm{\nabla u (x)} ^{ p - 2 } \nabla u ( x ) } = f
\end{equation}
for $ x \in \mathring{D} $ with
Dirichlet boundary conditions
where $ p \in (1, \infty) $
and where $ D = \{ x \in \R^d \colon \norm{ x } \leq 1 \} $
is the $ d $-dimensional unit ball.
It can be shown that the solution
of \cref{eq:nonlinear_DRM}
is a minimizer of the functional
$
  v \mapsto \frac{ 1 }{ p }
  \int_D ( \norm{ ( \nabla v )( x ) }^p - f v(x) ) \, \d x
$
over the space
$ W_0^{ 1, p }( \mathring{D} ) $
of Sobolev functions
with zero boundary values
and, therefore,
the \DRM\ is applicable to compute
an approximate of the solution
of \cref{eq:nonlinear_DRM}.
The boundary condition is again
incorporated using a penalty method.
In our simulations we use the values
$ d = 4 $, $ p = 9 $, and $ f = 1 $.
Note that the exact solution
$ u \colon D \to \R $
of
\cref{eq:nonlinear_DRM}
satisfies that
for all $ x \in D $
it holds that
$
  u(x)
  =
  p^{ - 1 } ( p - 1 )
  d^{ - 1 / ( p - 1 ) }
  ( 1 - \norm{x}^{ p / (p-1) } )
$.
In Figure~\ref{LaplacePowerRitz_[4, 32, 32, 32, 32, 1]_174133504536387}
we employ fully-connected feedforward \ANNs\ with
the \GELU\ activation and $ 4 $ hidden layers consisting of $ 32 $ neurons each
to approximate the exact solution
of \cref{eq:nonlinear_DRM}.
For the training we employ
mini-batches of size $ \num{256} $
and constant learning
rates of size $ \num{0.0003} $.
%
%
%
Furthemore, in
Figure~\ref{LaplacePowerRitz_[4, 32, 32, 32, 32, 1]_174133504536387}
we approximate the relative $ L^2( D; \R ) $-error
through a Monte Carlo approximation
with $ \num{16000} $ Monte Carlo samples
and we approximate
the $ L^1 $-error with respect
to the probability space
through a Monte Carlo approximation
with $ \num{50} $ independent simulations.
\begin{figure}[H]\caption{Deep Ritz for p-Laplace equation
	}\label{LaplacePowerRitz_[4, 32, 32, 32, 32, 1]_174133504536387}
	\includegraphics[width=\linewidth]{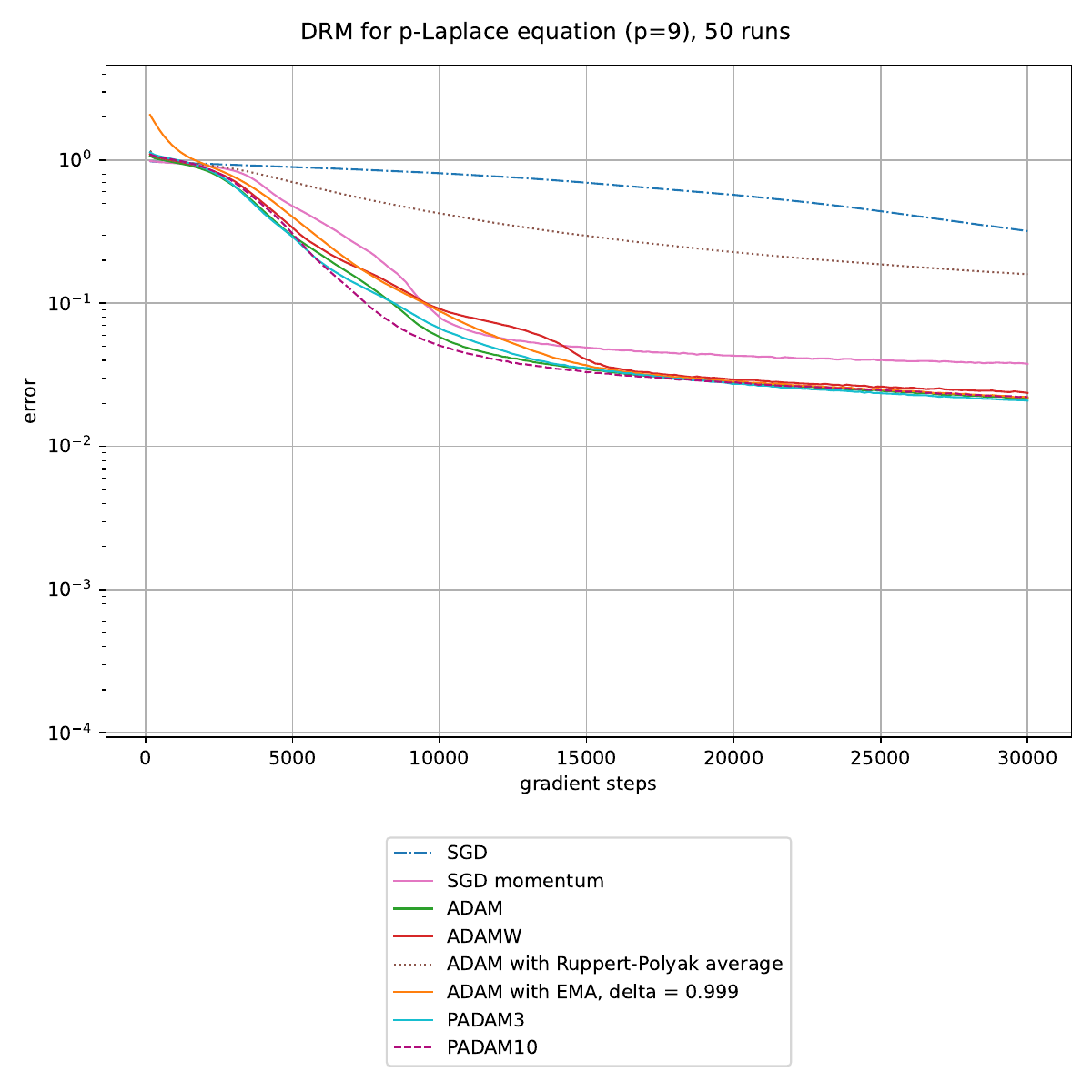}
\end{figure}

\subsection{Deep learning approximations
for optimal control problem}
\label{subsec:OC}

In the next example we consider
a controlled diffusion process
of the form
\begin{equation}
\label{eq:control_sde}
  \mathrm d X_t
  = -2 u_t \, \mathrm d t + \sqrt{2} \, \mathrm d W_t
\end{equation}
where $ W \colon [0,\infty) \times \Omega \to \R^d $
is a standard Brownian motion,
where $ X \colon [0, \infty ) \times \Omega \to \R^d $
is the diffusion process,
and where
$ u \colon [0, \infty ) \times \Omega \to \R^d $
is the control process.
We introduce the cost functional
\begin{equation}
\label{eq:cost_functional}
  J( t, x, u )
=
  \mathbb E \left[
    \int_t^T \|{u_s }\| ^2 \, \mathrm d s
    +
    \phi( X_T ) \, \middle \vert \, X_t = x
  \right] ,
\end{equation}
where the terminal cost function
$ \phi \colon \R^d \to \R $
satisfies for all $ x \in \R^d $ that
\begin{equation}
\label{eq:terminal_cost}
  \phi( x )
  = \ln\!\left(
    \tfrac{ 1 }{ 2 } ( \norm{x}^2 + 1 )
  \right)
  .
\end{equation}
We define the value function
$
  V(t, x ) = \inf_u J(t, x, u )
$
and attempt to compute the value
$
  V( 0, 0 )
$.
To approximate the solution
of the \SDE\ in \cref{eq:control_sde}
we consider a time discretization
of the form
$
  0 = t_0 < t_1 < \ldots < t_N = T
$
with
$
  \forall \, i \in \{0, 1, \ldots, N \} \colon t_i
  = \frac{iT}{N}
$.
The solution of \eqref{eq:control_sde}
is then approximated using a forward Euler method.
{Abbreviating $X_{t_n}$ by $X_n$ and $u_{t_n}$ by $u_n$}
for each $ n \in \{ 0, 1, \dots, N-1 \} $
we consider a control $ u_n $
of the form
$
  u_n \approx \mathcal N^{ \theta_n }( X_n )
$
where
$ \mathcal N^{ \theta_n } \colon \R^d \to \R^d $
is the realization function of
an \ANN\ with parameter vector $ \theta_n $.
Specifically, we use the values
$ d = 10 $, $ T = 1 $, $ N = 50 $,
and \ANNs\ with the \GELU\ activation
and $ 2 $ hidden layers consisting of $ 20 $ neurons each.
Additionally,
we employ batch normalization after the input layer and each hidden layer.
To approximately minimize the function 
$
u\mapsto J(0, 0, u ),
$
we use the time discretization and control
$
u_n \approx \mathcal N^{ \theta_n }( X_n )
$
introduced above and approximate the expectation in \cref{eq:cost_functional} through a Monte-Carlo approximation with $100000$ independently generated Brownian motion sample paths. 
For the training we use
mini-batches of size $ \num{256} $
and constant learning rates of size
$ \num{0.003} $.

We note that the value function
of the optimal control problem
in \cref{eq:control_sde},
\cref{eq:cost_functional},
and \cref{eq:terminal_cost}
satisfies
the \HJB\ equation
\begin{equation}
\label{eq:HJB_PDE_optimal_control}
\textstyle
  \frac{ \partial }{\partial t } V
  =
  \norm{\nabla_x V } ^2 - \Delta_x V ,
\qquad
  V( T, x ) = \phi( x )
\end{equation}
(cf., \eg, \cite[Chapters~3 and 4]{MR2533355}).
Moreover, we observe that
the Cole-Hopf transform
ensures that the solution of \cref{eq:HJB_PDE_optimal_control}
satisfies that
for all $ t \in [0,T] $,
$ x \in \R^d $ it holds that
\begin{equation}
\label{eq:cole-hopf}
  V( t , x )
  = -
  \ln\bigl(
    \E\bigl[
      \exp(
        - \phi( x + \sqrt{2} W_{ T - t } )
      )
    \bigr]
  \bigr),
\end{equation}
(cf., for example, \cite[Lemma 4.2]{EHanJentzen2017}).
To compute the error in Figure \ref{Optimal Control Problem_[10, 20, 20, 10]_174167700112340}
we employ \cref{eq:cole-hopf} to
approximately compute the value $ V( 0, 0 ) $
through a Monte Carlo approximation
with $ \num{40000000} $ Monte Carlo samples 
and display the relative errror $|(V_N-V(0,0))/V(0,0)|$ where 
$V_N\approx J(0,0,u)$ is the neural network approximation of $V(0,0)$ computed above.

Moreover, in Figure \ref{Optimal Control Problem_[10, 20, 20, 10]_174167700112340}
we approximate the $ L^1 $-error with respect to the probability space
through a Monte Carlo approximation with $ \num{10} $ independent simulations.

\begin{figure}[H]\caption{Optimal Control Problem
}\label{Optimal Control Problem_[10, 20, 20, 10]_174167700112340}
\includegraphics[width=\linewidth]{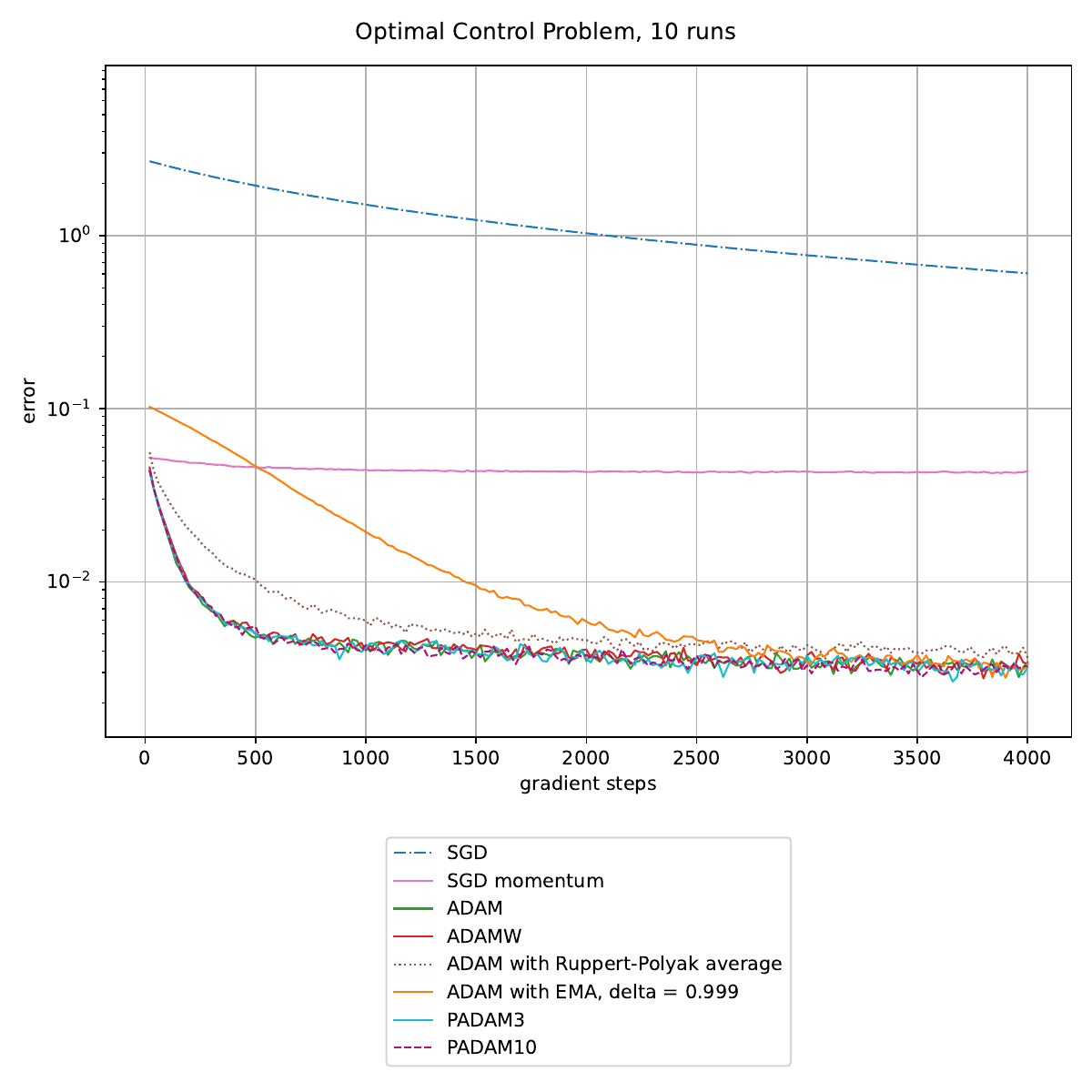}
\end{figure}

\subsection{Deep BSDE method
for Hamiltonian--Jacobi--Bellman equation}
\label{subsec:BSDE}

In the next example
we use the \deepBSDE\ method
introduced in
E et al.~\cite{EHanJentzen2017,HanJentzenE2018}
to approximately compute
the solution $ u \colon [0,T] \times \R^d \to \R $
of the \HJB\ equation
\begin{equation}
\label{eq:HJB}
  \tfrac{ \partial }{ \partial t } u
  + \Delta_x u
  = \norm{ \nabla_x u }^2
  ,
  \qquad
  u( T, x )
  = \ln\!\left( \tfrac{ 1 }{ 2 } ( \norm{x}^2 + 1 ) \right)
\end{equation}
for $ t \in [0,T] $, $ x \in \R^d $
on the spatial domain $ (-1, 1)^d $
at the initial time
$ ( u(0, x) )_{ x \in (-1,1)^d } $
where $ T = \nicefrac{ 1 }{ 5 } $
and where $ d = 10 $.
For the \deepBSDE\ method
we use a temporal discretization with $ 20 $ timesteps
and we approximate the gradient of the solution $u(t_n,x)$
of \cref{eq:HJB} at each time step $t_n=\tfrac{nT}{20}$ for $n=0,\dotsc,N-1$ 
by $\nabla_x u(t_n,x)\approx\mathcal N^{\theta_n}(x)$ where $\mathcal N^{\theta_n}\colon \R^d\to \R$ is the 
realization function of an \ANN\ with the \GELU\ activation 
and two hidden layers consisting of $30$ neurons each.
For the training we employ
mini-batches of size $ \num{256} $
and constant learning rates
of size $ \num{0.001} $.
To approximately compute the
{relative $ L^2( (-1,1)^d ; \R ) $-error of our neural network approximation $u(0,x)\approx\mathcal N^{\theta_0}(x)$} in
Figure~\ref{BSDEgen_[10, 30, 30, 10]_174220121537477}
we compare the computed results
with a reference solution computed
through a Monte Carlo approximation with
819200 Monte Carlo samples 
applied to
the Cole--Hopf transform
as in \cref{eq:cole-hopf}
(where we take the average over $400$ independent Monte Carlo approximation of \cref{eq:cole-hopf} using $2048$ Monte Carlo samples each).
In
Figure~\ref{BSDEgen_[10, 30, 30, 10]_174220121537477}
we approximate
the {relative
$ L^2( (-1,1)^d ; \R ) $-error}
through a
Monte Carlo approximation
with {\bf 1024} Monte Carlo samples
and
the $ L^1 $-error with respect to
the probability space is approximated
through a Monte Carlo approximation
with $ 50 $ independent simulations.

\begin{figure}[H]\caption{BSDE
}\label{BSDEgen_[10, 30, 30, 10]_174220121537477}
\includegraphics[width=\linewidth]{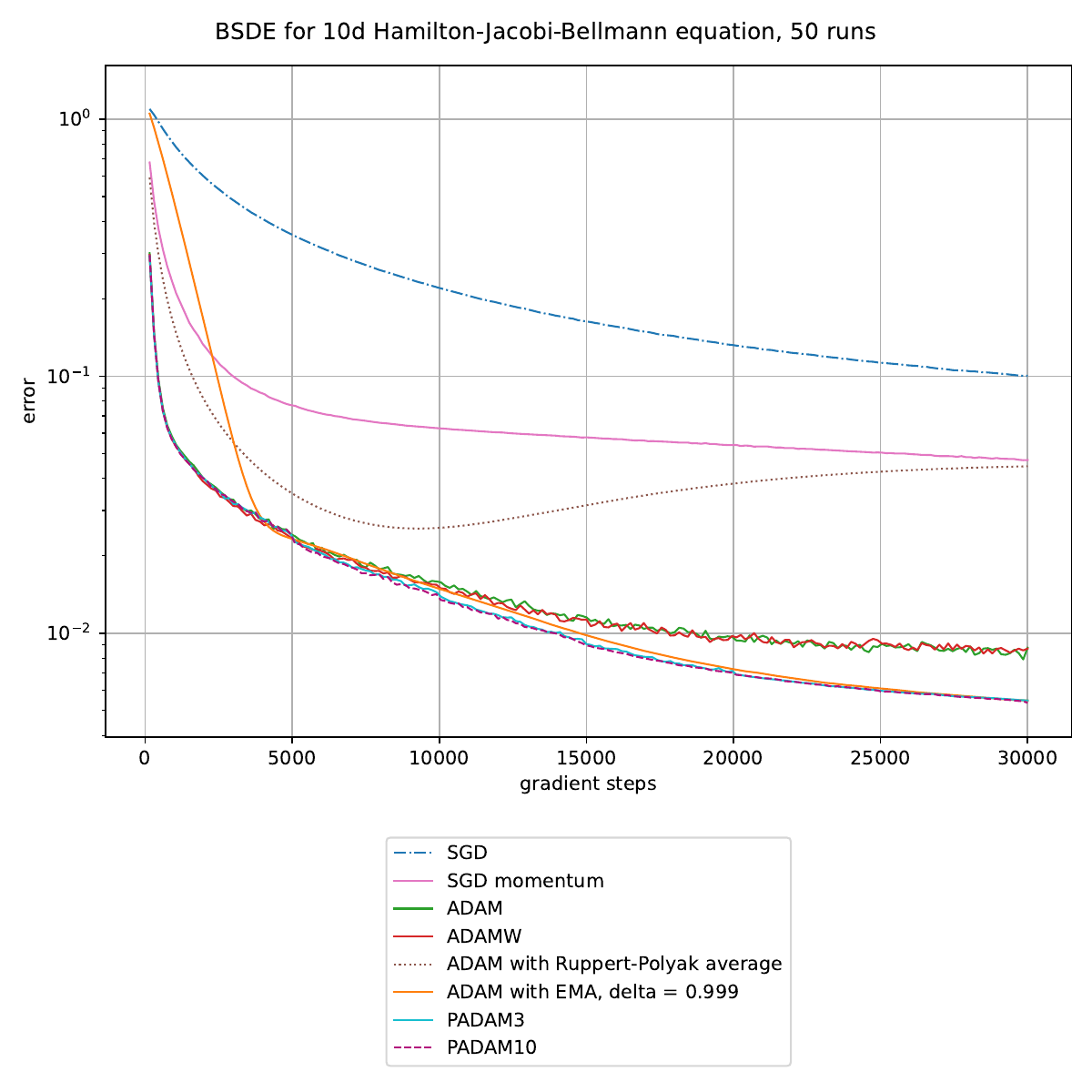}
\end{figure}
%
%

\subsection{Physics-informed neural networks (PINNs)
for Burgers equation}
\label{subsec:PINN1}

In the next example we employ the \PINN\ method
to approximately solve
the Burgers equation
\begin{equation}
\label{eq:pinn-burgers}
  \tfrac{ \partial }{ \partial t } u
  = \alpha \Delta_x u
  -
  u \,
  ( \tfrac{ \partial }{ \partial x } u ) ,
\qquad
  u(0, x )
  =
  \tfrac{
    2 \alpha \pi \sin( \pi x )
  }{
    \beta + \cos ( \pi x )
  }
\end{equation}
for $ t \in [0,T] $,
$ x \in \mathring{D} $
with Dirichlet boundary conditions
on the set $ D = [0,2] $
with the time horizon
$ T = \nicefrac{ 1 }{ 2 } $ and
where $ \alpha = \frac{ 1 }{ 20 } $
and $ \beta = \frac{ 11 }{ 10 } $.
Note that the exact solution
$ u \colon [0,T] \times D \to \R $
of \cref{eq:pinn-burgers} satisfies
that for all $ t \in [0,T] $, $ x \in D $
it holds that
$
  u(t,x)
  =
  \frac{
    2 \alpha \pi
    \sin( \pi x )
  }{
    \beta \exp( \alpha t \pi^2 )
    + \cos ( \pi x )
  }
$.
In
Figure~\ref{BurgersPINN1d_[2, 16, 32, 16, 1]_174109636670200}
we employ fully connected feedforward
\ANNs\ with the \GELU\ activation
and 3 hidden layers
consisting of $ 16 $, $ 32 $, and $ 16 $
neurons, respectively,
to approximate the solution
of \cref{eq:pinn-burgers}.
For the training
we employ mini-batches of size $ 256 $.
Moreover, for
\SGD\ and \SGD\ with momentum
we employ constant learning rates
of size $ \num{0.001} $
and for all other optimization methods
in
Figure~\ref{BurgersPINN1d_[2, 16, 32, 16, 1]_174109636670200}
we employ constant learning rates
of size $ \num{0.01} $.
In
Figure~\ref{BurgersPINN1d_[2, 16, 32, 16, 1]_174109636670200}
we approximate
the relative $ L^2( D; \R ) $-error
through a Monte Carlo approximation
with $ 10^5 $ Monte Carlo samples
and we approximate the $ L^1 $-error
with respect to the probability space
through a Monte Carlo approximation
with $ 50 $ independent simulations.

\begin{figure}[H]\caption{\PINN\ for Burgers equation
}\label{BurgersPINN1d_[2, 16, 32, 16, 1]_174109636670200}
\includegraphics[width=\linewidth]{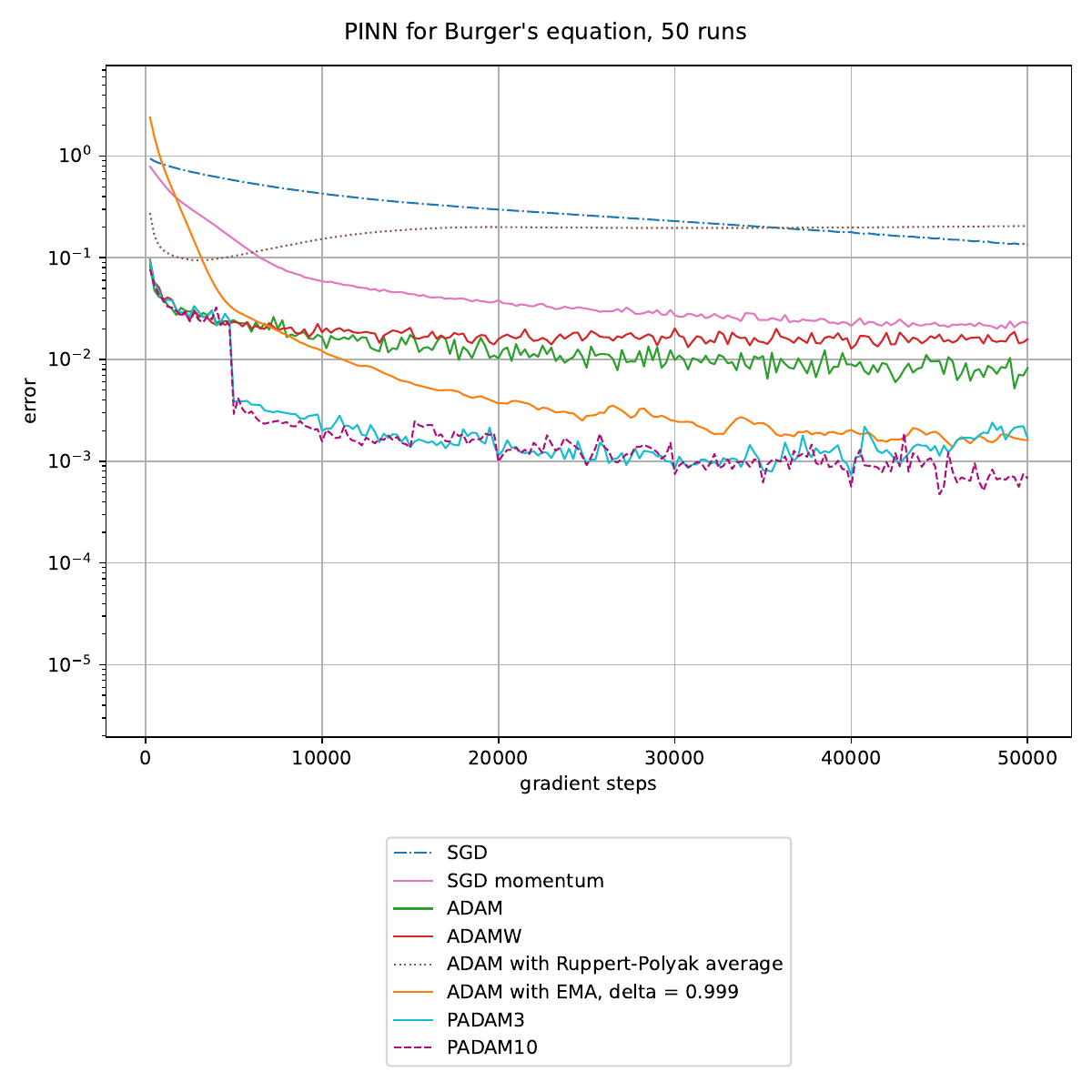}
\end{figure}
%
%
%

\subsection{PINNs for Allen--Cahn equation}
\label{subsec:PINN2}

In the next example we
employ the \PINN\ method
to aim to approximate
the solution
$
  u \colon [0, T] \times D \to \R
$
of the initial value Allen--Cahn \PDE\ problem
\begin{equation}
\label{eq:pinn-ac}
  \tfrac{ \partial }{ \partial t } u
  = \tfrac{1}{100} \Delta_x u + (u - u^3)
  ,
\qquad
  u(0, x )
  = \sin( \pi x_1 ) \sin( \pi x_2 )
\end{equation}
for $ t \in [0,T] $,
$ x = ( x_1, x_2 ) \in \mathring{D} $
with Dirichlet boundary conditions
on the set $ D = [0,2] \times [0,1] $
with the time horizon $ T = 4 $.
In
Figure~\ref{Allen Cahn PINN_[3, 32, 64, 32, 1]_174300226350547}
we employ fully connected feedforward
\ANNs\ with the \GELU\ activation
and $ 3 $ hidden layers consisting
of $ 32 $, $ 64 $, and $ 32 $ neurons,
respectively, to approximate
the solution of \cref{eq:pinn-ac}.
For the training we employ
mini-batches of size $ \num{256} $.
Moreover,
for \SGD\
we employ constant learning rates
of size $ \num{0.0001} $
(to avoid divergence),
for \SGD\ with momentum
we employ constant learning rates
of size $ \num{0.001} $
(to avoid divergence),
and for all other optimization methods
in Figure~\ref{Allen Cahn PINN_[3, 32, 64, 32, 1]_174300226350547}
we employ constant
learning rates of size $ \num{0.01} $.
To approximately compute the error
in Figure~\ref{Allen Cahn PINN_[3, 32, 64, 32, 1]_174300226350547}
we compare the computed results
with a reference solution computed by a
finite element method using $ 101^2 $
degrees of freedom in the spatial variable
and
$ \num{500} $ second order linear
implicit Runge-Kutta time steps.
In Figure~\ref{Allen Cahn PINN_[3, 32, 64, 32, 1]_174300226350547} the relative
$ L^2( D; \R ) $-error is approximated
through a Monte Carlo approximation
with $ \num{1000} $
Monte Carlo samples and
the $ L^1 $-error with respect
to the probability space is approximated
through a Monte Carlo approximation
with $ \num{50} $ independent simulations.

\begin{figure}[H]\caption{\PINN\ for Allen Cahn equation
	}\label{Allen Cahn PINN_[3, 32, 64, 32, 1]_174300226350547}
	\includegraphics[width=\linewidth]{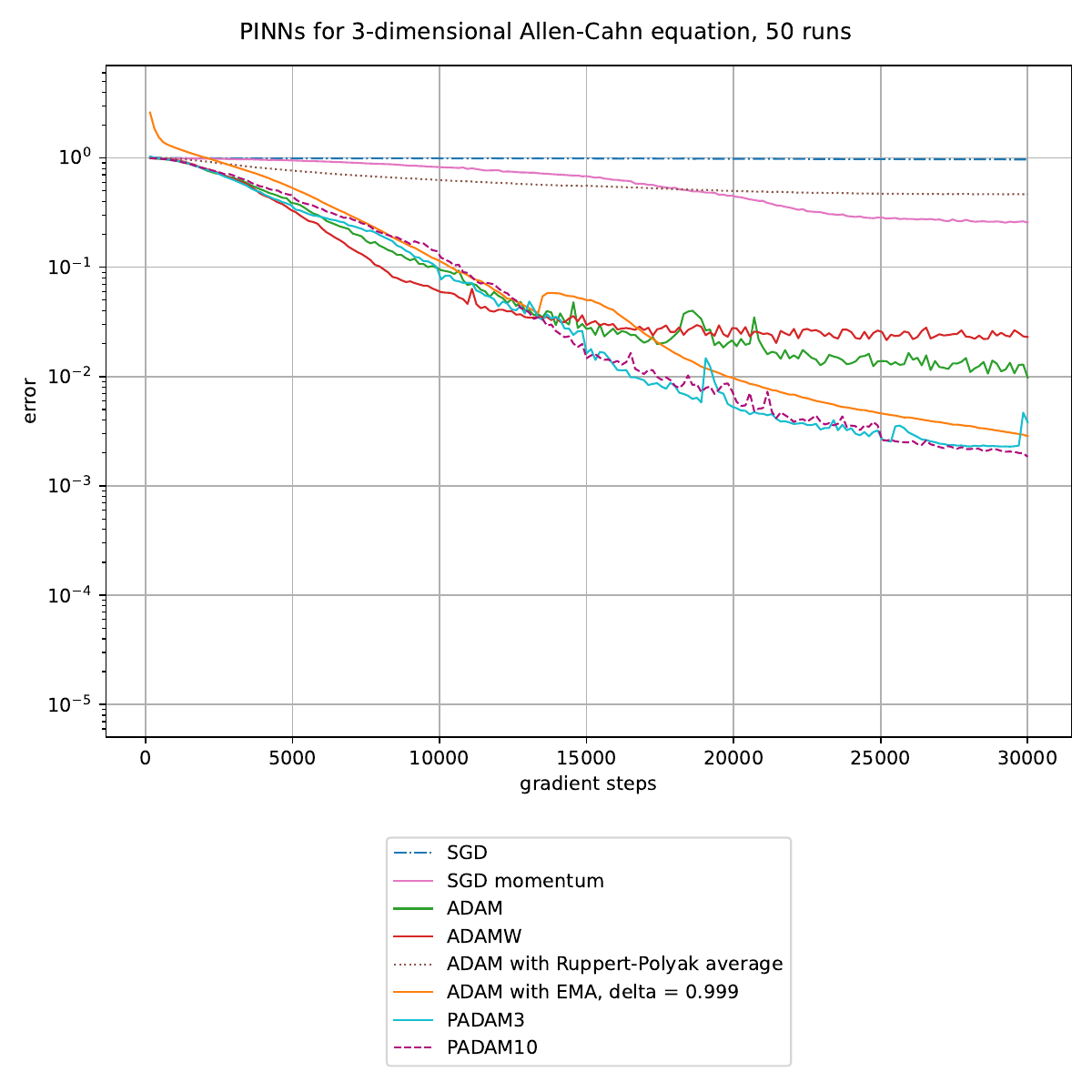}
	\end{figure}

\subsection{PINNs for Darcy flow}
\label{subsec:PINN3}

In the next example we apply the \PINN\ method
to approximately solve the parabolic linear Dacry flow \PDE
\begin{equation}
\label{eq:Darcy}
  - \operatorname{div}\bigl(
    a ( x ) ( \nabla u )( x )
  \bigr)
  = f
\end{equation}
for $ x \in \mathring{D} $
with Dirichlet boundary conditions
on the set $ D = [ -1 , 1 ]^2 $.
We consider $ f = 1 $
and
$ a \colon D \to \R $
satisfying
for all $ x = ( x_1, x_2 ) \in D $
that
$ a(x) = x_1 + 2 x_2 + 4 $.
In
Figure~\ref{Darcy PINN_[2, 32, 64, 32, 1]_174300226254314}
we employ fully connected feedforward \ANNs\ with
the \GELU\ activation
and 3 hidden layers consisting of
$ \num{32} $, $ \num{64} $,
and $ \num{32} $ neurons, respectively,
to approximate the solution of \cref{eq:Darcy}.
To compute the test error we compare the output with a reference solution obtained via a finite element method
using $ \num{16641} $ spatial degrees of freedom.
In
Figure~\ref{Darcy PINN_[2, 32, 64, 32, 1]_174300226254314}
we employ mini-batches of size $ \num{256} $.
Moreover,
for \SGD\ and \SGD\ with momentun
we employ constant learning rates
of size $ \num{0.0003} $
(to avoid divergence)
and for all other optimization methods
we employ constant learning rates
of size $ \num{0.003} $.
In Figure~\ref{Darcy PINN_[2, 32, 64, 32, 1]_174300226254314}
we approximate the relative
$ L^2( D; \R ) $-error
through a Monte Carlo approximation with
$ \num{3000} $
Monte Carlo samples and
we approximate the $ L^1 $-error
with respect to the probability space
through a Monte Carlo approximation
with $ \num{50} $ independent simulations.

\begin{figure}[H]\caption{PINNs for Darcy Flow
	}\label{Darcy PINN_[2, 32, 64, 32, 1]_174300226254314}
	\includegraphics[width=\linewidth]{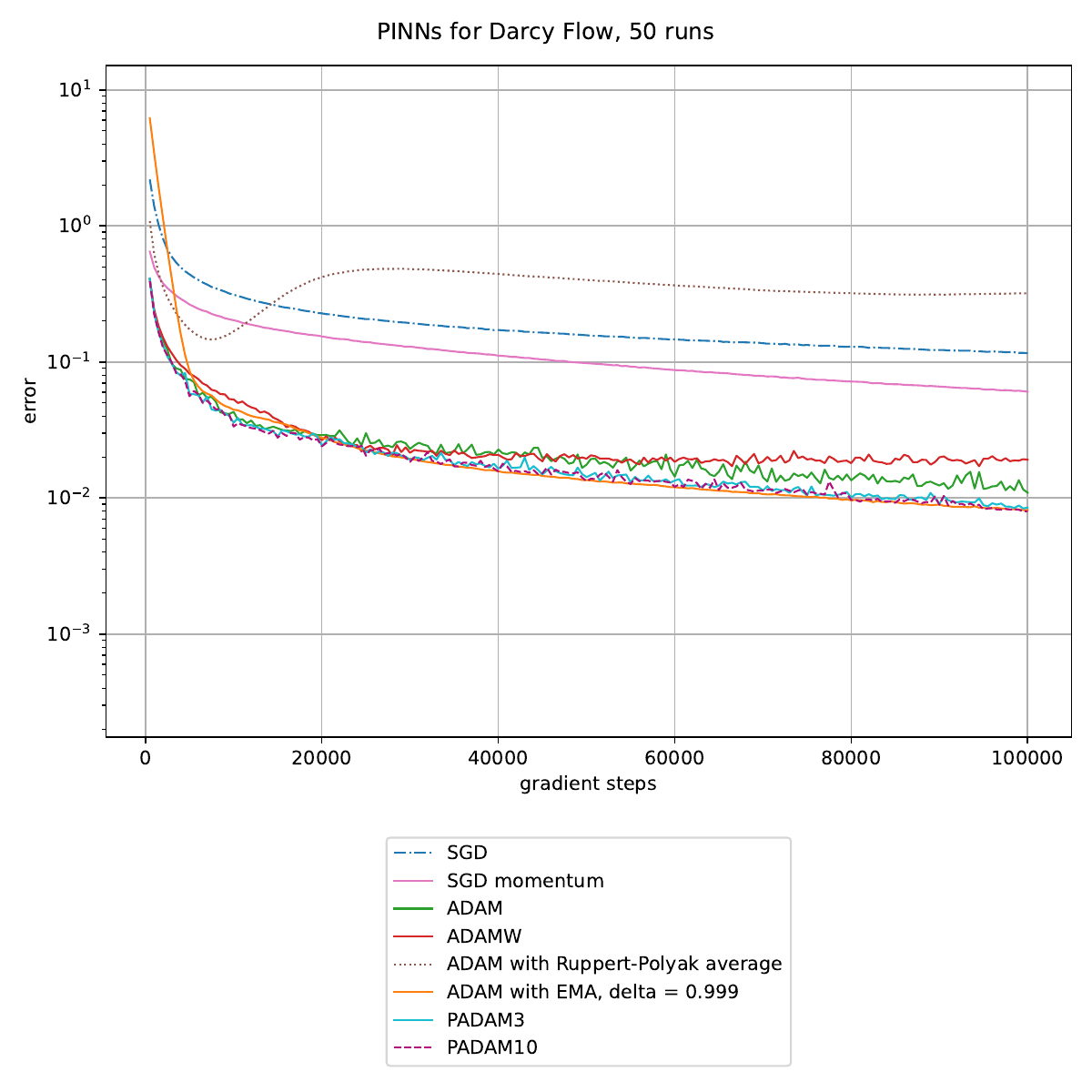}
\end{figure}

\subsection{Deep optimal stopping (DOS) method for American option}
\label{subsec:OS}

In the final example we employ
the \DOS\ method introduced
in Becker et
al.~\cite{MR4253974}
to approximately solve a $ 40 $-dimensional
optimal stopping problem
(cf.~\cite[Section 4.3.2.1]{MR4253974}).
Concretely, let $ d = 40 $,
$ r, \rho \in \R $,
$
  \beta = ( \beta_1, \dots, \beta_d )
$,
$
  \delta = ( \delta_1, \dots, \delta_d )
$,
$
  \xi = ( \xi_1, \dots, \xi_d ) \in \R^d
$
satisfy
for all $ i \in \cu{ 1, 2, \dots, d } $
that
$ r = \num{0.6} $,
$
  \rho = \frac{ 1 }{ d } \norm{ \beta }^2
$,
\begin{equation}
\textstyle
  \beta_i
  =
  \min\cu{
    \num{0.04} ( i - 1), \num{1.6} - \num{0.04} ( i - 1 )
  }
  ,
\quad
  \delta_i
  =
  r -
  \frac{ \rho }{ d }
  \rbr*{
    i - \frac{ 1 }{ 2 }
  }
  - \frac{ 1 }{ 5 \sqrt{d} }
  ,
\quad
  \xi_i = 100^{ 1 / \sqrt{d} } ,
\end{equation}
and
let $ \mu \colon \R^d \to \R^d $
and $ \sigma \colon \R^d \to \R^{ d \times d } $
satisfy for all $ x = ( x_1, \dots, x_d ) \in \R^d $
that
$
  \mu( x )
  =
  (
    ( r - \delta_1 ) x_1,
    \dots,
    ( r - \delta_d ) x_d
  )
$
and
$
  \sigma(x) =
  \operatorname{diag}
  ( \beta_1 x_1, \dots, \beta_d x_d )
$.
Let $ T = 1 $,
let
$
  ( \Omega, \cF, ( \cF_t )_{ t \in [0, T] }, \P )
$
be a filtered probability space,
let
$
  W = (W_t )_{ t \in [0,T] }
  \colon [0,T] \times \Omega \to \R^d
$
be a standard
$ ( \cF_t )_{ t \in [0,T] } $-Brownian motion,
and let
$
  X = (X_t)_{ t \in [0,T] }
  \colon [0, T ] \times \Omega \to \R^d
$
be an $ ( \cF_t )_{ t \in [0,T] } $-adapted
solution of the \SDE
\begin{equation}
  \d X_t =
  \mu ( X_t ) \, \d t + \sigma ( X_t ) \, \d W_t ,
  \qquad
  X_0 = \xi .
\end{equation}
Finally, let $ K = \num{95} $ and
define the payoff function
\begin{equation}
  g( s, x )
  =
  \exp( - r s )
  \textstyle
  \max\bigl\{
    K
    -
    \prod_{ k = 1 }^d \abs{ x_k }^{ 1 / \sqrt{d} }
    , 0
  \bigr\}
  .
\end{equation}
We are interested in approximating
the quantity
\begin{equation}
  \sup\cu[\big]{
    \E\br*{
      g( \tau, X_{ \tau } )
    }
    \in \R
    \colon
    \tau \colon \Omega \to [0,T] \text{ is an }
    ( \cF_t )_{ t \in [0,T] } \text{-stopping time}
  } ,
\end{equation}
which can be viewed as the price of an
American geometric-average put type option.
For this, we use the method described
in \cite{MR4253974}
with $ N = \num{100} $ time steps
and \ANNs\ with the \GELU\ activation and
two hidden layers consisting
of $ \num{240} $ neurons each.
Additionally, we employ batch normalization
after the input layer.
We use mini-batches of size
$ \num{256} $ and
constant learning rates of size
$ \num{0.0002} $.
To approximately compute
the error
in Figure~\ref{Stopping Net_[40, 240, 240, 1]_174133504483643},
we compare the result of the
different optimization methods with
the value $ \num{6.545} $,
which is calculated using a
binomial tree method
with \num{20000} nodes.
Furthemore,
in Figure~\ref{Stopping Net_[40, 240, 240, 1]_174133504483643}
the $ L^1 $-error with respect
to the probability space
is approximated through a
Monte Carlo approximation
with $ \num{10} $ independent simulations.

\begin{figure}[H]\caption{Optimal Stopping Problem
}\label{Stopping Net_[40, 240, 240, 1]_174133504483643}
\includegraphics[width=\linewidth]{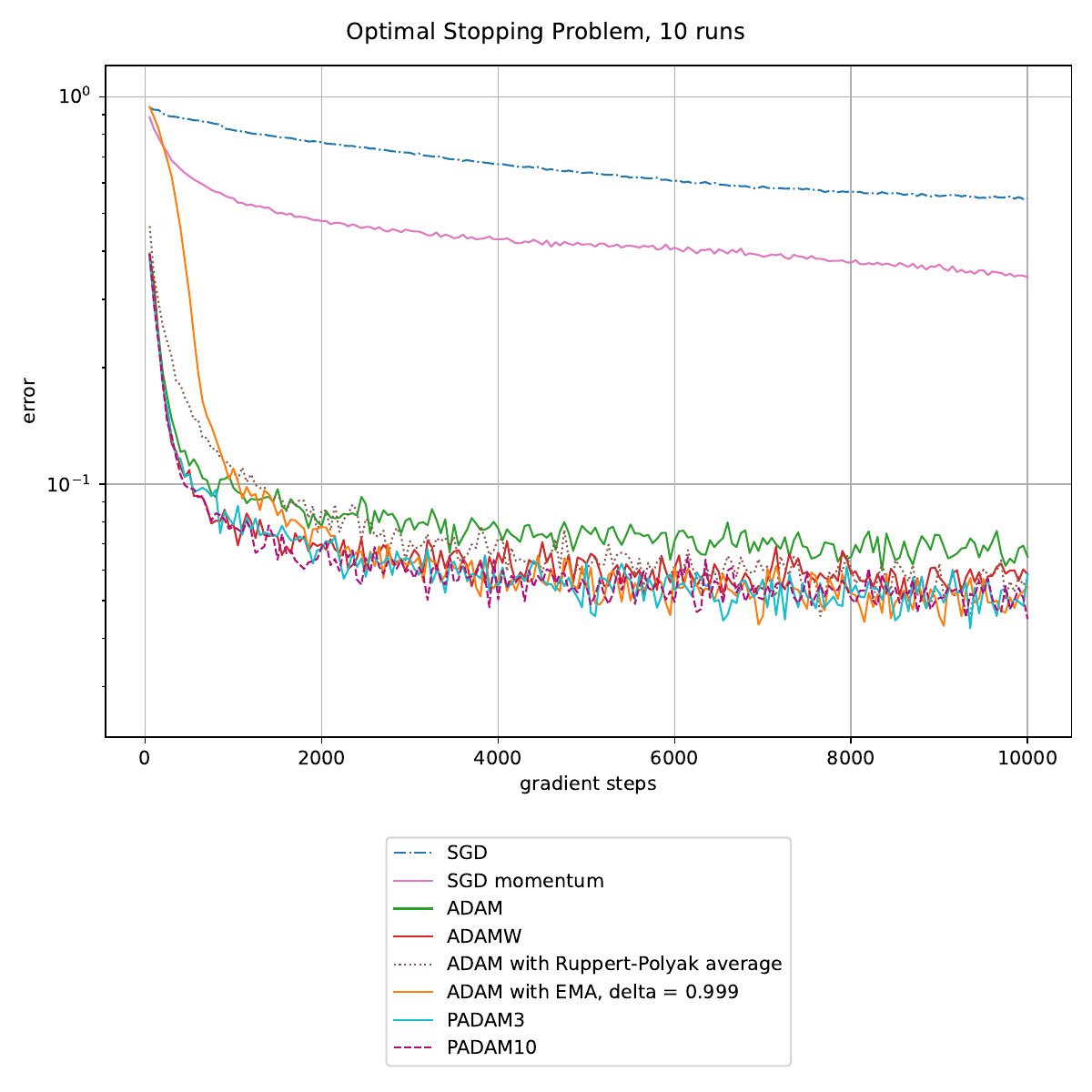}
\end{figure}

\section{Conclusion}
\label{sec:conclusion}

In this work we apply the proposed \Padam\ approach to a selection
of \numproblems stochastic optimization and deep \ANN\ learning problems
and compare it with some popular optimizers from the literature
such as standard \SGD, momentum \SGD, \Adam\ with
and without \EMA, and \AdamW.
In nearly all of the considered examples
\Padam\ achieves the smallest optimization error,
sometimes among others and sometimes exclusively.
%
%
We thus strongly suggest to consider \Padam\ and
related adaptive averaging techniques in the context of scientific machine learning problems.
In particular, this work aims to motivate further research for
suitable averaging procedures
when approximately solving \PDE, \OC, and related scientific computing problems
by means of deep learning methods.

\subsubsection*{Acknowledgments}

This work has been supported by the Ministry of Culture and Science NRW as part of the Lamarr Fellow Network.
In addition, this work has been partially funded by the Deutsche Forschungsgemeinschaft (DFG, German Research Foundation)
under Germany's Excellence Strategy EXC 2044-390685587, Mathematics Münster: Dynamics-Geometry-Structure.
Moreover, this work is supported via the AI4Forest project, which is funded by the German Federal Ministry of Education
and Research (BMBF; grant number 01IS23025A) and the French National Research Agency (ANR).
We also gratefully acknowledge the substantial computational resources that were made available to us by the
PALMA II cluster at the University of M\"{u}nster (subsidized by the DFG; INST 211/667-1).

\bibliographystyle{acm}
\bibliography{bibfile}



\end{document}